\documentclass[a4paper,12pt, reqno]{amsart}
\usepackage{amssymb,amsmath,latexsym}
\usepackage{color}
\usepackage{comment}
\allowdisplaybreaks

\newcommand{\bk}{\color{black}}

\def\1{{\bf 1}}
\def\ind{{\bf 1}}
\def\nn{\nonumber}
\def\bee{\begin{equation}}
\def\eee{\end{equation}}
 \def\sB {{\mathcal B}}

\def\R {{\mathbb R}}

\newtheorem{thm}{Theorem}[section]
\newtheorem{lemma}[thm]{Lemma}

\newtheorem{prop}[thm]{Proposition}
\newtheorem{corollary}[thm]{Corollary}
\newtheorem{remark}[thm]{Remark}

\numberwithin{equation}{section}

\def\qed{{\hfill $\Box$ \bigskip}}
\def\NN{{\mathcal N}}

\def\BB{{\mathcal B}}
\def\CC{{\mathcal C}}

\def\DD{{\mathcal D}}
\def\HH{{\mathcal H}}
\def\FF{{\mathcal F}}
\def\EE{{\mathcal E}}
\def\QQ{{\mathcal Q}}

\def\R{{\mathbb R}}

\def\E{{\mathbb E}}

\def\P{{\mathbb P}}

\def\wt{\widetilde}
\def\pf{\noindent{\bf Proof.} }

\usepackage{hyperref}
\begin{document}
\title[Markov processes with 
degenerate jump kernels]{
Harnack inequality and interior regularity for Markov processes with 
degenerate jump kernels}
\author{ Panki Kim \quad Renming Song \quad and \quad Zoran Vondra\v{c}ek}
\thanks{Panki Kim: This work was  supported by the National Research Foundation of
Korea(NRF) grant funded by the Korea government(MSIP) (No. NRF-2021R1A4A1027378)
}
\thanks{Renming Song: Research supported in part by a grant from
the Simons Foundation (\#429343, Renming Song)}
\thanks{Zoran Vondra\v{c}ek: Research supported in part by the Croatian Science Foundation under the project 4197.}

 \date{}

\begin{abstract}
In this paper we study interior potential-theoretic properties of purely discontinuous Markov processes in proper open subsets $D\subset \mathbb{R}^d$.The jump kernels of the processes may be degenerate at the boundary in the sense that they may vanish or blow up at the boundary. Under certain natural conditions on the jump kernel, we show that the scale invariant Harnack inequality holds for any proper open subset $D\subset \mathbb{R}^d$ and prove some interior regularity of harmonic functions. We also prove a Dynkin-type formula and several other  interior results.
\end{abstract}
\maketitle

\bigskip
\noindent {\bf  MSC 2020 Mathematical Sciences Classification System \bk}: Primary 60J45; Secondary 60J46, 60J76.



\section{Introduction and setting}\label{s:intro}

The goal of this paper is to
study interior potential-theoretic properties 
of purely discontinuous symmetric Markov 
processes in proper open subsets $D\subset \R^d$, $d\ge 1$.  
The main assumption is that we allow the jump 
kernels of the processes 
to be degenerate at the boundary $\partial D$. This includes 
the case when the  jump
kernels decay to zero at the boundary, as well as the case when they explode at the boundary. 
Examples of the former 
are subordinate killed L\'evy processes in smooth open sets $D$ studied in \cite{KSV19, KSV20}. An abstract approach to
 jump  kernels that decay 
 at the boundary is given in \cite[Section 3]{KSV21a}.
Compared with previous works, the main novelty of this paper
is that we also allow the possibility that the jump 
kernels blow up at the boundary. 
An example of such a case is the trace (or path-censored) process in $D$ of a nice isotropic L\'evy process in $\R^d$. In case $D$ is the half-space or an exterior $C^{1,1}$-open set, it can be deduced from
  \cite[Theorems 6.1 and 2.6]{BGPR} 
 that if $J(x,y)$ denotes the jump kernel of the trace process, then 
$\lim_{D\ni x\to z}J(x,y)=+\infty$ 
for all $z\in \partial D$ and $y\in D$. We will explain this example in much more detail in 
Subsection \ref{ss:E1} 
in the context of 
resurrected processes. 
A comprehensive study of potential-theoretic properties of 
such processes in the half-space with a  scale-invariant assumption is  given in \cite{KSV22b}, 
while the connection between these processes and positive self-similar Markov processes is given in \cite{KSV22a}.

We now describe our setup more precisely. 
Let $j:(0, \infty)\to (0, \infty]$ be a Borel function such that 
\begin{align}
\label{e:int}
\int_0^{\infty}\min(1,  r^2)j(r) r^{d-1}  dr<\infty.
\end{align}
We associate to $j$ an isotropic pure jump L\'evy process $X=(X_t, \P_x)$ in $\R^d$ with L\'evy measure $j(|x|)dx$. We further assume that 
\begin{equation}\label{e:assumption-on-j}
 j(r)\asymp r^{-d}\Psi(r)^{-1}, \quad \text{for all } r>0,
\end{equation}
where $\Psi$ is an increasing function satisfying  the following weak scaling condition: 
There exist constants $0<\delta_1\le \delta_2<1$ and $a_1, a_2>0$ such that 
\begin{equation}\label{e:H-infty}
   a_1(R/r)^{2\delta_1}\leq\frac{\Psi(R)}{\Psi(r)}\leq a_2(R/r)^{2\delta_2},\quad 
0<r<R <\infty.
\end{equation}
Here and throughout the paper, the notation $f\asymp g$ for 
non-negative functions
$f$ and $g$ means that there exists a constant $c\ge 1$ such that $c^{-1}g\le f \le c g$. 
A prototype of such a process $X$ is the isotropic $\alpha$-stable process in which case $\Psi(r)=r^{\alpha}$. This particular case already contains all the essential features of our results.

Throughout the paper, 
$a\wedge b:=\min \{a, b\}$, $a\vee b:=\max\{a, b\}$, and we use
$\delta_D(x)$ to denote the distance 
between $x$ and the boundary $\partial D$. 
Let, also, $j(x,y):=j(|x-y|)$.

For a given proper open subset $D$ of $\R^d$,
we will consider a process  on $D$
associated with a pure jump Dirichlet form whose jump kernel has the form 
\begin{align}
\label{e:DJ}
J(x,y)=\sB(x,y)j(|x-y|).
\end{align}
Here $\sB(x,y)$ may depend on $\delta_D(x), \delta_D(y)$ and $|x-y|$, and is allowed to 
vanish at the boundary or to explode at the boundary. 
The main concern of this paper is on interior results, so we do not
need to impose any regularity assumption on the boundary $\partial D$.

We remark here that, 
when $D$ is bounded, 
it suffices to assume that the function $\Psi$ satisfies
the property \eqref{e:H-infty} 
for all $0 < r \le R  \le \mathrm{diam}(D)$.
Indeed, in this case, $\Psi$ can be extended to satisfies
the property \eqref{e:H-infty} 
for all $0 < r \le R  <\infty$ trivially. Then for any Borel function $j:(0, \mathrm{diam}(D)]\to (0, \infty]$  such that 
$ j(r)\asymp r^{-d}\Psi(r)^{-1}$ for all $ 0<r \le \mathrm{diam}(D)$,  we extend $j$ so that 
\eqref{e:int} and \eqref{e:assumption-on-j} hold.

Throughout this paper, we will assume that 
$\sB: D \times D \to [0, \infty)$  satisfies the following hypothesis:

\medskip
\noindent
\textbf{(H1)} $\sB(x,y)=\sB(y,x)$ for all $x,y\in D$.

\medskip

\medskip
\noindent
\textbf{(H2)} For any $a\in (0,1)$ there exists $C_1=C_1(a)\ge 1$ such that for all $x,y\in D$ satisfying
 $\delta_D(x)\wedge \delta_D(y)\ge a|x-y|$,
 it holds that
$$
C_1^{-1}\le \sB(x,y)\le C_1.
$$ 
Without loss of generality, we assume that $a\mapsto C_1(a)$ is decreasing on $(0,1)$.

\medskip
\noindent
\textbf{(H3)} For any $a >0$ there exists 
$C_2=C_2(a)>0$
such that
\begin{align}
\label{e:H30}
\int_{D, |y-x|>a\delta_D(x)} J(x,y)dy \le 
C_2 \Psi(\delta_D(x))^{-1}\, .
\end{align}
The assumption  \textbf{(H3)}
states that the tail of the jump measure depends only on the distance to the boundary of $D$  and $\Psi$ (or $j$), 
and clearly provides 
sufficient integrability of 
the function $y \mapsto J(x, y)$ away from the point $x$. 
 We also note that it follows from \textbf{(H2)} that $D\ni x\mapsto \sB(x, x)$ is bounded between two positive constants.

The assumption  \textbf{(H3)} clearly holds when $\sB(x,y)$ is 
bounded above by a positive constant,  see \eqref{e:HH}. 
In Subsection \ref{ss:E2}, we give examples of  $\sB(x,y)$, satisfying \textbf{(H3)}, that may explode at the boundary.

Assumptions \textbf{(H2)}-\textbf{(H3)} are scale invariant. For some of the results their weaker and non-scale invariant versions will suffice. Therefore we introduce 

\medskip
\noindent
\textbf{(H2-w)} 
For any relatively compact open subset $U$ of $D$, there exists a constant 
$C_{3}=C_{3}(U)\ge 1$ 
such that $C_{3}^{-1}\le \sB(x,y)\le C_{3}$ for all $x,y\in U$.

\medskip
\noindent
\textbf{(H3-w)} For any relatively compact open set $U\subset D$ and open set $V$ with $\overline{U}\subset V\subset D$,
\begin{align}\label{e:H3}
\sup_{x\in U}\int_{D\setminus V}J(x,y)\, dy  
<\infty.
\end{align}

It is easy to deduce that \textbf{(H2)}, respectively \textbf{(H3)}, implies \textbf{(H2-w)}, respectively \textbf{(H3-w)}, see Lemma \ref{l:local-boundedness}.

For functions $u,v:D\to \R$, define
\begin{equation}\label{e:df-uv}
\EE^D(u,v):=\frac12 \int_D\int_D (u(x)-u(y))(v(x)-v(y))J(x,y)\, dy\, dx.
\end{equation} 
Assumptions \textbf{(H1)}, \textbf{(H2-w)}-\textbf{(H3-w)} are sufficient to conclude, 
see Lemma \ref{l:df-well-defined}, that $\EE^D(u,v)$ is well defined for all $u,v\in C_c^{\infty}(D)$.
By Fatou's lemma, 
$( \EE^D, C_c^{\infty}(D))$ 
is closable in $L^2(D, dx)$. Let $\FF^D$ be the closure of $C_c^{\infty}(D)$ under $\EE^D_1:=\EE^D+(\cdot, \cdot)_{L^2(D,dx)}$. 
Then $(\EE^D, \FF^D)$ is a regular Dirichlet form on $L^2(D, dx)$. 

Let $\kappa:D\to [0,\infty)$ be a function satisfying
\begin{equation}\label{e:bound-on-kappa}
\kappa(x)\le C_4 \frac{1}{\Psi(\delta_D(x))}\, ,\quad x\in D\, ,
\end{equation}
for some constant $C_4>0$.
Then $\kappa$ is locally bounded in $D$.
Set
$$
\EE^{D, \kappa}(u,v):=\EE^D(u,v)+\int_D u(x)v(x)\kappa(x)\, dx\, .
$$
Since $\kappa$ is locally bounded, the measure $\kappa(x)dx$ is a positive Radon measure  charging no set of zero capacity. Let $\FF^{D,\kappa}:=\wt{\FF^D}\cap L^2(D, \kappa(x)dx)$, where $\wt{\FF^D}$ is the family of all quasi-continuous functions on $\FF^D$. By \cite[Theorems 6.1.1 and 6.1.2]{FOT}, 
$(\EE^{D, \kappa},\FF^{D,\kappa})$ 
is a regular Dirichlet form on $L^2(D, dx)$ having $C_c^{\infty}(D)$ as a special standard core.
Let $((Y_t^{\kappa})_{t\ge0}, (\P_x)_{x\in D\setminus \NN})$
be the associated Hunt 
process with lifetime $\zeta$, where $\NN$ is an exceptional set.
We add a cemetery point $\partial$ to the state space $D$ and define $Y_t^{\kappa}=\partial$
for $t\ge \zeta$. We will write $D_{\partial}=D\cup \{\partial\}$. 
Any function $f$ on $D$ is automatically extended to $D_{\partial}$
by setting $f(\partial)=0$.
In Section \ref{s:regularization},
we will show that we can remove the exceptional set $\NN$ so the process $Y^{\kappa}$ can start 
from every point in $D$, see Proposition \ref{p:YDregularized}.

Our process may not be Feller but the next hypothesis will allow us to establish a Dynkin-type formula on any relatively compact open set of $D$, see Theorem \ref{t:Dynkin}.

\medskip
\noindent
\textbf{(H4)}  
If $\delta_2 \ge 1/2$, then there exists $\theta>2\delta_2-1$ 
with the property 
that for any $a>0$ there exists 
$C_5=C_5(a)>0$  such that
$$
|\sB(x, x)-\sB(x,y)|\le C_5\left(\frac{|x-y|}{\delta_D(x)\wedge \delta_D(y)}\right)^{\theta}\,
\quad \text{ for all }x,y\in D \text{ with }
\delta_D(x)\wedge \delta_D(y)  \ge 
 a  |x-y|.
$$

The final hypothesis ensures that jumping from two points close to each other to a faraway point is comparable:

\medskip
\noindent
\textbf{(H5)} For any  $\epsilon \in (0,1)$ there exists $C_6=C_6(\epsilon)\ge 1$ with the following property: For all $x_0\in D$ and $r>0$ with $B(x_0, (1+\epsilon)r)\subset D$, we have
$$
C_6^{-1}\sB(x_1,z)\le \sB(x_2,z) \le C_6 \sB(x_1, z)\, ,\quad \text{for all }x_1,x_2\in B(x_0,r), \, \,z\in D\setminus B(x_0, (1+\epsilon)r)\, .
$$

An immediate consequence of \textbf{(H5)} is the following: For any
$\epsilon \in (0,1)$ there exists $C_7=C_7(\epsilon)\ge 1$ with the property that,  for all $x_0\in D$ and $r>0$ with $B(x_0, (1+\epsilon)r)\subset D$, we have
\begin{equation}\label{e:H5-for-J}
C_7^{-1}J(x_1,z)\le J(x_2,z) \le C_7 J(x_1, z)\, ,\quad \text{for all }x_1,x_2\in B(x_0,r), \, \,z\in D\setminus B(x_0, (1+\epsilon)r)\, .
\end{equation} 
Indeed, by \textbf{(H5)} 
we have $ J(z,x_1)=\sB(z,x_1)j(|z-x_1|)\le 
C_6 \sB(z,x_2) j(|z-x_1|)$. Since  $|x_1-z|\le |x_2-z|+|x_1-x_2|\le |x_2-z|+2r\le |x_2-z| +(2/\epsilon)|x_2-z|=(1+2/\epsilon) |x_2-z|$, it follows from  \eqref{e:H-infty} that $j(|z-x_1|)\le c_1 j(|z-x_2|)$. This proves 
\eqref{e:H5-for-J}.

We now
compare hypotheses \textbf{(H1)}-\textbf{(H5)} with the 
hypotheses 
\cite[\textbf{(B1)}-\textbf{(B5)}]{KSV21a}. 
The symmetry hypothesis \textbf{(H1)} is the same as
\cite[ \textbf{(B1)}]{KSV21a}. 
Also, \textbf{(H4)}-\textbf{(H5)} are precisely 
 \cite[\textbf{(B4)}-\textbf{(B5)}]{KSV21a}. 
 The key difference is that here we do not assume that $\sB(x,y)$ is bounded from above on $D\times D$ by a
positive constant which was 
\cite[\textbf{(B2)}]{KSV21a}. 
Instead, we assume \textbf{(H2)} which is a two-sided version of 
\cite[\textbf{(B3)}]{KSV21a}. 
 Finally, \textbf{(H3)}, which implies that $J(x,y)$ is sufficiently integrable, is automatically satisfied under 
 the boundedness condition \cite[\textbf{(B2)}]{KSV21a}.

Recall that a Borel function $f$
defined on $D$ is said to be \emph{harmonic} in an open set $U\subset D$ with respect to the process $Y^{\kappa}$ if, for every bounded open set $V\subset \overline{V}\subset U$, it holds that
$$
\E_x[|f(Y^{\kappa}_{\tau_V})|; \tau_V<\infty]<\infty \quad \text{and} \quad  f(x)=\E_x[f(Y^{\kappa}_{\tau_V}); \tau_V<\infty], \quad \text{for all } x\in V,
$$
where $\tau_V=\inf\{t>0:\, Y^{\kappa}_t \notin V\}$ is the first exit time from $V$.

Here are our main results under the scale invariant hypotheses \textbf{(H1)}-\textbf{(H5)}. The first one is the scale invariant 
Harnack inequality.

\begin{thm}[scale invariant Harnack inequality]\label{t:uhp*}
Suppose $D$ is 
a proper open subset of $\R^d$ and assume that 
\textbf{(H1)}-\textbf{(H5)},
 \eqref{e:assumption-on-j}-\eqref{e:H-infty} and \eqref{e:bound-on-kappa} hold. 

\begin{itemize}
\item[(a)] 
There exists a constant $C_8>0$ such that 
for any $r\in (0,1]$, $B(x_0, r) \subset D$ and any non-negative function $f$ in $D$ which is 
harmonic in $B(x_0, r)$ with respect to $Y^{\kappa}$, we have
$$
f(x)\le C_8 f(y), \qquad \text{ for all } x, y\in 
B(x_0, r/2).
$$

\item[(b)] 
There exists a constant $C_9>0$ 
such that for any  $L>0$, any $r \in (0, 1]$, 
all $x_1,x_2 \in D$ with $|x_1-x_2|<Lr$ and $B(x_1,r)\cup B(x_2,r) \subset D$ 
and any non-negative function $f$ in $D$ which is  harmonic in $B(x_1,r)\cup B(x_2,r)$ with respect to $Y^{\kappa}$ we have
$$
f(x_2)\le C_9 C_1 
(\tfrac{1}{2(L+1)})
L^{d+2\delta_2} f(x_1)\, .
$$
\end{itemize}
\end{thm}

The second result is H\"older continuity of bounded 
  harmonic functions.

\begin{thm}\label{t:holder}
Suppose $D$ is a proper open subset of $\R^d$
and assume that 
\textbf{(H1)}-\textbf{(H4)}, 
\eqref{e:assumption-on-j}-\eqref{e:H-infty} and \eqref{e:bound-on-kappa} hold. 
Then there exist $C_{10}>0$  and $\beta>0$ 
such that for any
$(r, x_0)\in (0, \infty)\times D$ 
with $B(x_0, 6r)\subset D$ and
any  bounded
 function in $D$ which is harmonic in $B(x_0, r)$ with respect to $Y^{\kappa}$, 
$$
|f(x)-f(y)|\le C_{10}\|f\|_\infty \left(\frac{|x-y|}{r}\right)^\beta, 
\qquad  \text{for all } x, y\in B(x_0, r/2).
$$
\end{thm}
Study of Harnack inequality and regularity of harmonic functions of general discontinuous processes started with the paper \cite{BL}. There have been many papers on this subject since then, see \cite{BK, CKW19, CKW16b, CKW16a, DK, K2, KW1,KW2, SV04} and the references therein.
Almost all these papers assume that the jump kernels
of the processes are non-degenerate.  
Some of the exceptions are \cite{KSV19, KSV20, KSV21a, KSV21b},  where the jump kernels are allowed to decay to zero at the boundary. 
 The main assumption of this paper is \textbf{(H3)}, 
which, together with \textbf{(H2)}, implies 
an upper bound on the tail of the jump measure for $r <\delta_D(x)$.
 A similar condition appeared in \cite[Definition 1.5]{CKL} in the study of law of iterated logarithm for general Markov processes, see also \cite{CKW19} for  a global version of this condition. In our setup, the jump kernel is allowed to be 
degenerate at the boundary of $D$, including the possibility of decay to 0 or blow up to infinity. 
The case of the jump kernel decaying to 0 at the boundary has been studied in \cite{KSV19, KSV20, KSV21a, KSV21b}.  
The possibility of the jump kernel blowing up at the boundary 
leads to some complications in proving  the results above. 

 We note that results obtained in this paper, in particular Theorem \ref{t:uhp*}, will be used \cite{KSV22b}.

Organization of the paper: 
In Section \ref{s:prelim-res} we collect some preliminary results that follow from \textbf{(H1)} and \textbf{(H2)}-\textbf{(H3)},
respectively \textbf{(H2-w)}-\textbf{(H3-w)},
 which allow the construction of the process.

In Section \ref{s:regularization} we first show that, for any relatively compact open set $U\subset D$, the Dirichlet forms 
of the killed processes $X^U$ and $Y^{\kappa, U}$ are comparable. 
Using this and some pretty involved analysis, we then show
that $Y^{\kappa, U}$ can be identified with a process which can 
start from any point in $U$. This allows us to remove the
exceptional set $\NN$, see Proposition \ref{p:YDregularized}.

Section \ref{s:generator} is devoted to the study of the generator of the process $Y^{\kappa}$. In order to handle the singularity of the jump kernel, we need hypothesis \textbf{(H4)}, see the proof of Proposition \ref{p:operator-interpretation}. Different from all previous works, 
it turns out that, 
due to the possible blow-up of the jump kernel at the boundary, the action of the generator on a function compactly supported in $D$ need not be bounded. 
This makes the task of proving the Dynkin-type formula (Theorem \ref{t:Dynkin}) difficult.

In Section \ref{s:hi} we establish all necessary ingredients for the proof of Harnack's inequality including the exit time estimates from balls (Proposition \ref{p:exit-time-estimate})  and Krylov-Safonov-type estimate (Lemma \ref{l:A2}), 
and give the proofs of Theorems \ref{t:uhp*} and \ref{t:holder}. 
We also give a sketch of the proof of a non-scale invariant Harnack inequality (Proposition \ref{t:uhp}).

Building up on  standard theory and 
some results from Section \ref{s:regularization},  we show in Section \ref{s:exist-gf} 
that, in the transient case, $Y^{\kappa}$ has a Green function. We also 
prove the natural result
that if the killing function $\kappa$ is strictly positive, then $Y^{\kappa}$ is transient.

Finally, in Section \ref{s:examples}, we give several families of jump kernels satisfying our hypotheses. The main examples are trace processes and more general resurrection processes given in Subsection \ref{ss:E1}. 
In Subsection \ref{ss:E2}, we first give 
examples  of jump 
kernels satisfying our hypotheses, which may blow up at the boundary, 
and then we look at the setting of 
\cite[Section 4]{KSV22b} in case of the  half-space,  and show that all hypotheses are satisfied. 
This section
can be read independently of the rest of the paper and  provides further motivation for studying jump kernels exploding at the boundary through some concrete examples. Some readers may want to  glance through it before reading the main body of the paper.

Throughout this paper, the positive constants 
$\delta_1$, $\delta_2$, $\theta$, $a_1$ and $a_2$
will remain the same.
We will use the following convention:
Capital letters $C_i, i=1,2,  \dots$ will denote constants
in the statements of results and assumptions. The labeling of these constants will remain the same. Lower case letters 
$c_i, i=1,2,  \dots$ are used to denote the constants in the proofs
and the labeling of these constants starts anew in each proof.
The notation $c_i=c_i(a,b,\ldots)$, $i=0,1,2,  \dots$ indicates  constants depending on $a, b, \ldots$.


\section{Preliminary results}\label{s:prelim-res}
Throughout this section we assume that 
\textbf{(H1)} and 
\eqref{e:assumption-on-j}-\eqref{e:H-infty} hold.
First note that by \eqref{e:assumption-on-j}-\eqref{e:H-infty},
\begin{align}
 \label{e:ww1}
\int_{|x-y|>a} j(x,y)dy
& \le c_1 \Psi(a)^{-1} a^{2\delta_1}
\int_{|x-y|>a}  |x-y|^{-d-2\delta_1}dy\le 
c_2\Psi(a)^{-1}
\end{align}
and 
\begin{align}
 \label{e:ww2}
\int_{|x-y|<a} |x-y|^2 j(x,y)dy
& \le c_3 \Psi(a)^{-1} a^{2\delta_2}
\int_{|x-y|<a}  |x-y|^{-d-2\delta_2+2}dy  \le 
c_4a^2\Psi(a)^{-1}.
\end{align}

We will use the following notation: 
For $U\subset D$, 
$d_U:=\mathrm{diam}(U)$ and $\delta_U:=\mathrm{dist}(U, \partial D)$.

\begin{lemma}\label{l:local-boundedness}
(a) If \textbf{(H2)} holds, then \textbf{(H2-w)} also holds with 
$C_3=C_1\left(\frac{\delta_U}{d_U+\delta_U}\right) \ge 1$.

\noindent (b) If \textbf{(H3)} holds, then \textbf{(H3-w)} also holds.
\end{lemma}
\pf (a) Let $U$ be a relatively compact open subset of $D$. For $x,y\in U$, we have 
$$
\delta_D(x)\wedge \delta_D(y)\ge \delta_U\ge \frac{\delta_U}{d_U}|x-y|\ge
 \left(\frac{\delta_U}{d_U+\delta_U} \right)|x-y|.
$$ 
Hence, with 
$a:=\frac{\delta_U}{d_U+\delta_U}\in (0,1)$, 
we get $C_1(a)^{-1}\le \sB(x,y)\le C_1(a)$ for all $x,y\in U$. 

\noindent (b) Let $U$ be a relatively compact open subset of $D$ and $V$ an open set with $\overline{U}\subset V\subset D$. Let 
$a:=\mathrm{dist}(U, D\setminus V)/\sup_{z\in U}\delta_D(z)$.
Then for all $x\in U$ and $y\in D\setminus V$, it holds that 
$|y-x|>\mathrm{dist}(U, D\setminus V)=a\sup_{z\in U}\delta_D(z)> a\delta_D(x)$.  
Therefore, by \textbf{(H3)}, we have that 
$$
\int_{D\setminus V}J(x,y)\, dy
\le \int_{D, |y-x|> a\delta_D(x)}J(x,y)\, dy \le 
c_1(a)
\Psi(\delta_D(x))^{-1}, \quad x\in U,
$$
implying that
$$
\sup_{x\in U}\int_{D\setminus V} J(x,y)\, dy \le 
c_1(a)
\Psi(\delta_U)^{-1}<\infty.
$$
Hence, \textbf{(H3-w)} holds.
\qed

\begin{lemma}\label{l:H3-strong} 
Suppose \textbf{(H2)}-\textbf{(H3)} hold. There exists a constant $C_{11}>0$ such that for all $x\in D$ and $r\in (0, \delta_D(x)]$,
\begin{equation}\label{e:kest}
\int_{D, \, |y-x|>r}J(x,y)\, dy \le C_{11}\Psi(r)^{-1}.
\end{equation}
\end{lemma}
\pf If $|y-x| \le \delta_D(x)/2$, then $\delta_D(y) \ge 
\delta_D(x)-|y-x| \ge \delta(x)/2$ so
$|y-x| \le \delta_D(y)\wedge \delta_D(x)$. 
Thus, by \textbf{(H2)}-\textbf{(H3)} and \eqref{e:ww1}, 
for  $r \le \delta_D(x)$,
\begin{align*}
\int_{D, |y-x|>r/2}J(x,y)\, dy &\le 
\int_{D, |y-x|>\delta_D(x)/2}J(x,y)dy+C_1
\int_{D, \delta_D(x)/2 \ge |y-x|>r/2}j(x,y)\, dy \nn\\
&\le C_2(1/2)\Psi(\delta_D(x))^{-1} +C_1
\int_{ |y-x|>r/2}j(x,y)\, dy
\nn\\  &\le C_2(1/2)\Psi(r)^{-1} +c_1
\Psi(r)^{-1}  =c_2 \Psi(r)^{-1} .
\end{align*}
\qed

We note here that if $\sB(x,y)\le c_1$ for all $x,y\in D$, then $J(x,y)\le c_1j(|y-x|)$, and thus 
by \eqref{e:ww1},
\begin{align}
& \int_{|y-x|>a\delta_D(x)} J(x,y)dy\le
c_1 \int_{|y-x|>a\delta_D(x)} j(|y-x|)dy\le 
c_2(a) \Psi(\delta_D(x))^{-1}.\label{e:HH}
\end{align}
Therefore, \textbf{(H3)} holds true. This fact was already mentioned in the introduction.

\begin{lemma}\label{l:consequences}
Suppose that  \textbf{(H2-w)} and \textbf{(H3-w)} hold.

\noindent
(a) For any relatively compact open subset $U$ of $D$,
$$
\sup_{x\in U}\int_D (1\wedge |x-y|^2)J(x,y)\, dy <\infty.
$$

\noindent
(b) For any compact set $K$ and open set $V$ with $K\subset V\subset D$, 
\begin{equation}\label{e:fot-df}
\iint_{K\times K}|x-y|^2 J(x,y)\, dy\, dx <\infty, \qquad \int_K \int_{D\setminus V}J(x,y)\, dy\, dx <\infty.
\end{equation}

\end{lemma}
\pf (a) Let $V$ be a relatively compact open set such that $\overline{U}\subset V\subset D$. 
By \textbf{(H3-w)},
we only need to check $\sup_{x \in U}\int_V (1\wedge |x-y|^2)J(x,y)\, dy<\infty.$

Since $V$ is relatively compact, 
by  \textbf{(H2-w)}, 
$\sB(x,y)\le c_1$ for $x, y \in V$ for some $c_1=c_1(V)$.
Therefore, using \eqref{e:ww2}, for $x \in U$,
\begin{align*}
\int_V (1\wedge |x-y|^2)J(x,y)\, dy &\le c_1 \int_V |x-y|^2 J(x,y)\, dy \le c_1 \int_V |x-y|^2 j(|x-y|)\, dy \\
& \le c_1 \int_{B(x, d_V)}|x-y|^2 j(|x-y|)\, dy \le  c_2 d_V^2\Psi(d_V)^{-1} <\infty.
\end{align*}

\noindent 
(b) Let $U$ be a relatively compact open set such that $K\subset U\subset \overline{U}\subset V\subset D$. For $x,y\in K$, it holds that $|x-y|^2\le (d_K\vee 1)^2(1\wedge |x-y|^2)$. Therefore,
\begin{align*}
& \int_K \int_K |x-y|^2 J(x,y)\, dy\, dx \le (d_K\vee 1)^2 \int_K \int_K (1\wedge|x-y|^2) J(x,y)\, dy\, dx\\
\quad & \le (d_K\vee 1)^2 \int_K \int_D (1\wedge|x-y|^2) J(x,y)\, dy\, dx \\
\quad & \le (d_K\vee 1)^2 |K| \sup_{x\in U} \int_D (1\wedge|x-y|^2) J(x,y)\, dy <\infty,
\end{align*}
where the finiteness of the integral follows from part (a). 
For the second integral in \eqref{e:fot-df},
let $b:=\mathrm{dist}(K, D\setminus V)$. Then for $x\in K$, $y\in D\setminus V$, $|x-y|\ge b\ge b\wedge 1$. Therefore,
\begin{align*}
& \int_K \int_{D\setminus V}J(x,y)\, dy\, dx \le \frac{1}{(b\wedge 1)^2}\int_K \int_{D\setminus V}(1\wedge |x-y|^2)J(x,y)\, dy\, dx\\
\quad & \le \frac{1}{(b\wedge 1)^2}\int_K \int_{D}(1\wedge |x-y|^2)J(x,y)\, dy\, dx\\
\quad & \le \frac{1}{(b\wedge 1)^2} |K| \sup_{x\in K}\int_{D}(1\wedge |x-y|^2)J(x,y)\, dy <\infty
\end{align*}
by part (a). 
\qed

Recall that $\EE^D$ is defined in \eqref{e:df-uv}.
Condition \eqref{e:fot-df} is  sufficient and necessary 
for $\EE^D(u,u)<\infty$ for all $u\in C_c^{\infty}(D)$, see \cite[p.7]{FOT}. Therefore, 
under \textbf{(H1)}, 
\textbf{(H2-w)}-\textbf{(H3-w)}, 
$\EE^D(u,u)$ is finite for all $u\in C_c^{\infty}(D)$. In particular, \eqref{e:df-uv} is well defined for all $u,v\in C_c^{\infty}(D)$.
In fact, we will need a little bit more.

\begin{lemma}\label{l:df-well-defined}
For all $u\in C_c^2(\R^d)$ and $v\in C_c^2(D)$, 
$$
\int_D \int_D |(u(x)-u(y))(v(x)-v(y))| J(x,y)\, dy\, dx <\infty.
$$
\end{lemma}
\pf
Let $K=\mathrm{supp}(v)$ and $V$ be a relatively compact open subset of $D$ with  $K \subset V\subset \overline V \subset D$. Then
\begin{align*}
& \int_D \int_D |(u(x)-u(y))(v(x)-v(y))| J(x,y)\, dy\, dx \\
\qquad & =\int_V \int_V +\int_{D\setminus V}\int_V +\int_V \int_{D\setminus V}+\int_{D\setminus V}\int_{D\setminus V}=:I+II+III+IV.
\end{align*}
By \eqref{e:fot-df}, we have
$$
I\le \|\nabla u\|_{\infty}\|\nabla v\|_{\infty}
\int_{\overline V}\int_{\overline V}
 |x-y|^2 J(x,y)\, dy\, dx <\infty.
$$
Next,
$$
II=\int_{D\setminus V}\int_K |(u(x)-u(y))v(y)|J(x,y)\, dy\, dx \le 2\|u\|_{\infty}\|v\|_{\infty}\int_{D\setminus V}\int_K J(x,y)\, dy\, dx<\infty
$$
again by \eqref{e:fot-df}. The integral $III$ is estimated in the same way as $II$, while $IV=0$. \qed


\section{Regularization of the process}\label{s:regularization}

In this section we will show that, under \textbf{(H1)}, \textbf{(H2-w)}, \textbf{(H3-w)}, \eqref{e:assumption-on-j}-\eqref{e:H-infty} and
the condition that for any relatively compact open set $U$,  
\begin{equation}\label{e:bound-on-kappa0}
\|\kappa |_U\|_\infty < \infty\, ,
\end{equation}
we can remove the exceptional set $\NN$ and so the process $Y^{\kappa}$ can start from every point $x\in D$.
For this purpose, we 
will use an auxiliary process $Z$ on $\R^d$, 
with jump kernel $J_\gamma$ defined below.
The process $Z$ can start from every point in $\R^d$.
We will first prove a result stating that, for a relatively compact open subset $U$ of $D$, 
the Dirichlet forms of the parts of the processes $X$ and $Y^{\kappa}$ on $U$ 
are comparable. Recall from 
Section \ref{s:intro}
that $X$ is a L\'evy process in $\R^d$ with L\'evy measure $j(|x|)dx$, so that 
its jump kernel  is precisely $j(x, y)$.

For a relatively compact open subset $U$ of $D$,
let $Y^{\kappa, U}$ be the process $Y^{\kappa}$ killed upon exiting $U$, that is, the part of the process $Y^{\kappa}$ in $U$. The Dirichlet form of $Y^{\kappa, U}$ is 
$(\EE^{D, \kappa}, \FF^{D,\kappa}_U)$, 
where $\FF^{D,\kappa}_U=\{u\in \FF^{D,\kappa}: \, u=0 \textrm{ q.e.~on  } D\setminus U\}$. 
Here q.e.~means that the equality holds quasi-everywhere, that is, except on a set of capacity zero with respect to $Y^{\kappa}$. 
Let \begin{equation}\label{e:kappaU}
\kappa^U(x)=\int_{D\setminus U} J(x,y)\, dy  \quad \text{ and } \quad \kappa_U(x)=\kappa^U(x) +\kappa(x)\, , \quad x\in U\, .
\end{equation}
Then, 
for $u,v\in \FF^{D,\kappa}_U$,
\begin{equation}\label{e:EDkappa}
\EE^{D, \kappa}(u,v)=\frac{1}{2}\int_U\int_U (u(x)-u(y))(v(x)-v(y))J(x,y)\, dy\, dx +
\int_U u(x)v(x)\kappa_U (x)\, dx\, .
\end{equation}
Note that it follows from \textbf{(H3-w)}
and \eqref{e:bound-on-kappa0}
that $\kappa_U(x)<\infty$ for all $x\in U$. 
Further, since $C_c^{\infty}(D)$ is a special standard core of $(\EE^{D, \kappa}, \FF^{D,\kappa})$,  $C_c^{\infty}(U)$ is a core of 
$(\EE^{D, \kappa}, \FF^{D,\kappa}_U)$.

For $u,v:\R^d\to \R$, let
\begin{eqnarray*}
\QQ(u,v)&:=&\frac{1}{2}\int_{\R^d} \int_{\R^d} (u(x)-u(y))(v(x)-v(y))j(|x-y|)\, dy \, dx\, ,\\
\DD(\QQ)&:=& \{u\in L^2(\R^d, dx): \QQ(u,u)<\infty\}.
\end{eqnarray*}
Then $(\QQ, \DD(\QQ))$ is the regular Dirichlet form corresponding to $X$. Let $X^U$ 
denote the part of the process $X$ in $U$. 
The Dirichlet form of $X^U$ is $(\QQ^U, \DD_U(\QQ))$, where
\begin{equation}\label{e:Qkappa}
\QQ^{U}(u,v)=\frac{1}{2}\int_U \int_U (u(x)-u(y))(v(x)-v(y))j(|x-y|)\, dy \, dx+\int_U  u(x)v(x) 
\kappa^X_U(x)\, dx, 
\end{equation}
\begin{equation}\label{e;kappaU-X}
\kappa^X_U(x)
:=\int_{ \R^d  \setminus U}j(|y-x|)\, dy, \qquad x\in U
\end{equation}
and $\DD_U(\QQ)=\{u\in \DD(\QQ) :\, u=0 \textrm{ q.e.~on  } \R^d\setminus U\}$.
Here q.e.~means, except on a set of capacity zero with respect to $X$.

Recall that $\delta_U=\mathrm{dist}(U, \partial D)$ and $d_U=\mathrm{diam}(U)$.
By \textbf{(H2-w)}, there exists a constant $c_1=c_1(U)\ge 1$ such that $c_1^{-1}\le \sB(x,y)\le c_1$ for all $x,y\in U$. 
We note that if the stronger \textbf{(H2)} holds, then by Lemma \ref{l:local-boundedness}(i),  
the constant $c_1$ is equal to $C_1\left(\frac{\delta_U}{d_U+\delta_U}\right)$, and thus only depends on $\frac{\delta_U}{d_U+\delta_U}$. 
This fact that will be important in Lemma \ref{l:df-comparison1}.
Together with \eqref{e:assumption-on-j}, 
the boundedness of $\sB(\cdot, \cdot)$ on $U\times U$
implies that there 
exist  $c_2>0$  and $c_3>0$ such that
\begin{equation}\label{e:JD-in-U}
\frac{c_2}{|x-y|^d\Psi(|x-y|)}\le J(x,y)\le \frac{c_3}{|x-y|^d\Psi(|x-y|)}\, ,\quad x,y\in 
U\, .
\end{equation} 
This can be written equivalently as 
\begin{equation}\label{e:equal-domains-1}
c_4^{-1} j(|x-y|)\le J(x,y)\le c_4 j(|x-y|)\, ,\quad x,y\in U\, ,
\end{equation}
for some $c_4\ge 1$.
Let $V$ be the  $\delta_U/2$\bk-neighborhood of $U$, 
that is, $V:=\{x\in D:\, \mathrm{dist}(x, U)<\delta_U/2\}$. Then
\begin{equation}\label{e:JD-in-U2}
\kappa_U(x)= \kappa^U(x)+\kappa(x)\
=\int_{D\setminus V}J(x,y)\, dy+\int_{V\setminus U}J(x,y)\, dy +\kappa(x),\, 
 \qquad x\in U.
\end{equation} 
Similarly as above we conclude that $c_5^{-1} j(|x-y|)\le J(x,y)\le c_5 j(|x-y|)$ for all $x,y\in V$ with 
$c_5:=c_5(U)\ge 1$.
If the stronger \textbf{(H2)} holds, 
then the constant is equal to $C_1\left(\frac{\delta_V}{d_V+\delta_V}\right)=C_1\left(\frac{\delta_U}{2d_U+5\delta_U}\right)\le  C_1\left(\frac{1}{5}\frac{\delta_U}{d_U+\delta_U}\right)$, and thus depends only on $\frac{\delta_U}{d_U+\delta_U}$. 
Moreover, by \textbf{(H3-w)}, $\sup_{x\in U}\int_{D\setminus V}J(x,y)\, dy =:c_6\ < \infty$, with 
$c_6=c_6(U)$.
By setting 
$c_{7}:=
\|\kappa_{|U}\|_\infty$, 
we get
$$
c_5^{-1}\int_{V\setminus U}j(|x-y|)\, dy \le \kappa_U(x)\le c_6 +c_5 \int_{V\setminus U}j(|x-y|)dy+c_{7}\, ,\quad x\in U\, .
$$
Since 
$$
\inf_{x\in U}\int_{V\setminus U}j(|x-y|)\, dy \ge |V\setminus U|j(\mathrm{diam}(V))=:c_{8}>0\, ,
$$
we conclude that
$$
c_5^{-1}\int_{V\setminus U}j(|x-y|)\, dy \le \kappa_U(x)\le c_{9} \int_{V\setminus U}j(|x-y|)\, dy\, ,\quad x\in U\, .
$$
Further, since
\begin{equation}\label{e:kappa-U-X-new}
\kappa_U^X(x)=\int_{\R^d\setminus V}  j(|x-y|)\, dy+\int_{V\setminus U} j(|x-y|)\, dy\, ,\quad x\in U\, 
\end{equation}
and $\sup_{x\in U}\int_{\R^d\setminus V}j(|x-y|)\, dy=:c_{10}<\infty$, we see that there is a constant $c_{11}>0$ such that 
$$
\int_{V\setminus U}j(|x-y|)\, dy \le 
\kappa_U^X(x)
\le c_{11} \int_{V\setminus U}j(|x-y|)\, dy\, ,\quad x\in U\, .
$$
It follows that 
\begin{equation}\label{e:kappaU-vs-kappaU-X}
c_9^{-1}\kappa_U(x)\le \kappa_U^X(x)\le c_{11}c_5\kappa_U(x),
\end{equation}
with constants $c_5$, $c_9$ and $c_{11}$ depending on $U$.

Let $\mathrm{Cap}^{Y^{\kappa,U}}$ and $\mathrm{Cap}^{X^U}$ denote the capacities with respect to the killed processes $Y^{\kappa,U}$, and $X^U$ respectively. 
\begin{lemma}\label{l:df-comparison}
Assume that \textbf{(H1)}, \textbf{(H2-w)}, \textbf{(H3-w)}, 
\eqref{e:assumption-on-j}-\eqref{e:H-infty} 
and \eqref{e:bound-on-kappa0} hold.
Let $U$ be a relatively compact open subset  of $D$.
\noindent
(a) There exists a constant 
$C_{12}=C_{12}(U)\ge 1$ such that
\begin{equation}\label{e:comp-dfs}
C_{12}^{-1}\EE^{D,\kappa}(u,u)\le \QQ(u,u)\le C_{12}\EE^{D,\kappa}(u,u)\qquad \textrm{for all }u\in C_c^{\infty}(U).
\end{equation}

\noindent (b)
For any Borel $A\subset U$, 
\begin{equation}\label{e:comp-cap}
C_{12}^{-1}\mathrm{Cap}^{Y^{\kappa,U}}(A)\le \mathrm{Cap}^{X^U}(A)\le C_{12} \mathrm{Cap}^{Y^{\kappa, U}}(A),
\end{equation}
where $C_{12}$ is the constant from part (a).
\end{lemma}
\pf (a) This follows immediately from \eqref{e:EDkappa}, \eqref{e:Qkappa}, \eqref{e:equal-domains-1} and \eqref{e:kappaU-vs-kappaU-X}.

\noindent (b) 
Since $ C_c^{\infty}(U)$ is a core for both $(\QQ^U, \DD_U(\QQ))$ and $(\EE^{D, \kappa}, \FF^{D,\kappa}_U)$ 
by using the definition of capacity as in \cite[2.1]{FOT}, the claim follows from part (a). 
\qed

\begin{lemma}\label{l:df-comparison1}
Assume that \textbf{(H2)}-\textbf{(H3)}, \eqref{e:assumption-on-j}-\eqref{e:H-infty} and \eqref{e:bound-on-kappa} hold.
Let $U$ be a relatively compact open subset  of $D$.
Then the constant  $C_{12}$ 
in Lemma \ref{l:df-comparison} depends only on 
$\frac{\delta_U}{d_U+\delta_U}$  and is decreasing in $\frac{\delta_U}{d_U+\delta_U}$. 
\end{lemma}
\pf 
Let $V$ be the  $\delta_U/2$\bk-neighborhood of $U$. 
Recall that by \textbf{(H2)}, 
Lemma \ref{l:local-boundedness}(a) and \eqref{e:assumption-on-j},
we have  with 
$c_1>1$, depending on $U$ only through $\frac{\delta_U}{d_U+\delta_U}$ and being a decreasing function of $\frac{\delta_U}{d_U+\delta_U}$,  
so that 
\begin{equation}\label{e:JD-in-V}
\frac{c_1^{-1}}{|x-y|^d\Psi(|x-y|)}\le J(x,y)\le 
\frac{c_1}{|x-y|^d\Psi(|x-y|)}\, ,\quad x,y\in V.
\end{equation} 
By Lemma \ref{l:H3-strong}, for all $x\in U$, 
$$
\int_{D\setminus V}J(x,y)\, dy  \le \int_{D,\,  |y-x|>\delta_U/2} J(x,y)\, dy\le 
\frac{C_{11}}
{\Psi(\delta_U)}.
$$
Using \eqref{e:bound-on-kappa} and the fact that 
$\Psi(\delta_D(x)) \ge \Psi(\delta_U)$ for $x \in U$, we get
$$
c_1^{-1}\int_{V\setminus U}j(|x-y|)\, dy \le \kappa_U(x)\le\frac{C_{11}}
{\Psi(\delta_U)} +c_1 \int_{V\setminus U}j(|x-y|)dy +  \frac{C_4}{\Psi(\delta_U)} \, ,\quad x\in U\, .
$$
If $x\in U$ and $y\in V$, then $|x-y|\le d_U+\delta_U$, hence $j(|x-y|)\ge j(d_U+\delta_U)$. Moreover, we can find a point $z$ so that 
$B(z,\delta_U/4)\subset V\setminus U$, 
to obtain that 
$$
\inf_{x\in U}\int_{V\setminus U}j(|x-y|)\, dy\ge |B(z, (\delta_U/4)|j(d_U+\delta_U)\ge 
c_2\frac{\delta_U^d }
{(d_U+\delta_U)^d \Psi(d_U+\delta_U)},
$$
where $c_2$ is independent of $x$ and $U$. 
Thus, 
$$
c_1^{-1}\int_{V\setminus U}j(|x-y|)\, dy \le \kappa_U(x)\le
(c_1 +c_3\frac
{(d_U+\delta_U)^d \Psi(d_U+\delta_U)}{\delta_U^d \Psi(\delta_U)})  \int_{V\setminus U}j(|x-y|)dy \, ,\quad x\in U\, ,
$$
with $c_3=(C_{11}+C_4)/c_2$.

Recall that
$$
\kappa_U^X(x)=\int_{\R^d\setminus V}  j(|x-y|)\, dy+\int_{V\setminus U} j(|x-y|)\, dy\, ,
\quad x\in U.
$$
Since $\R^d\setminus V \subset \R^d\setminus B(x, \delta_U/2)$ for $x\in U$,  
by \eqref{e:ww1} we have
$$
\int_{\R^d\setminus V}j(|x-y|)\, dy\le \int_{\R^d\setminus B(x, \delta_U/2)}j(|x-y|)\, dy \le 
 c_4 \Psi(\delta_U)^{-1}
$$
for some $c_4$ independent of $x$ and $U$. 
Thus 
$$
\int_{V\setminus U}j(|x-y|)\, dy \le 
\kappa_U^X(x)
\le (1 +c_5\frac
{(d_U+\delta_U)^d \Psi(d_U+\delta_U)}{\delta_U^d \Psi(\delta_U)}) \int_{V\setminus U}j(|x-y|)\, dy\, ,\quad x\in U\, ,
$$
with $c_5=c_4/c_2$. 
It follows that 
$$
(c_1 +c_3\frac
{(d_U+\delta_U)^d \Psi(d_U+\delta_U)}{\delta_U^d \Psi(\delta_U)}) ^{-1}\kappa_U(x)\le \kappa_U^X(x)\le c_{1}
(1 +c_5\frac
{(d_U+\delta_U)^d \Psi(d_U+\delta_U)}{\delta_U^d \Psi(\delta_U)}) 
\kappa_U(x),
$$
where $c_1$ depends on $U$ only through $\frac{\delta_U}{d_U+\delta_U}$ and is a decreasing function of $\frac{\delta_U}{d_U+\delta_U}$,  
and $c_3$ and $c_5$ are independent of $U$.
This 
and \eqref{e:H-infty} imply 
that 
\begin{equation}\label{e:kappaXX}
c_6^{-1}\kappa_U(x)\le \kappa_U^X(x)\le c_{6}\kappa_U(x),
\end{equation}
where $c_6$ depends on $U$ only through $\frac{\delta_U}{d_U+\delta_U}$ and is a decreasing function of $\frac{\delta_U}{d_U+\delta_U}$.
Using \eqref{e:EDkappa}, \eqref{e:Qkappa} and \eqref{e:JD-in-V}, 
the statement of the lemma follows from \eqref{e:kappaXX} in the same way as in the proof of Lemma \ref{l:df-comparison}. 
\qed

In the remainder of this
section we assume that \textbf{(H1)}, \textbf{(H2-w)}, \textbf{(H3-w)}, \eqref{e:assumption-on-j}-\eqref{e:H-infty} and \eqref{e:bound-on-kappa0} hold.

\begin{lemma}\label{l:YDkappaU}
Let $U$ be a relatively compact 
open subset of $D$. 
The process $Y^{ \kappa, U}$ can be refined to start from every point in $U$. Moreover, it is strongly Feller.
\end{lemma}
\pf 
Define a kernel $J_\gamma(x,y)$ on $\R^d\times \R^d$ by $J_\gamma(x,y)=J(x,y)$ for $x,y\in U$, and $J_{\gamma}(x,y)=\gamma  j(|x-y|)$ otherwise, where $\gamma>0$ is a positive constant 
to be chosen later.
Using $J_\gamma$, we define
$$
\CC(u,u):=\frac12 \int_{\R^d}\int_{\R^d} (u(x)-u(y))^2 
J_\gamma(x,y)\, dx\, dy \text{ and }
\DD(\CC):=\{u\in L^2(\R^d):\, \CC(u,u)<\infty\}\, .
$$
Note that $C_c^{\infty}(\R^d)$ is a special standard core of $ \DD(\CC)$.
By \eqref{e:assumption-on-j} 
and \eqref{e:JD-in-U}, $
J_\gamma(x,y) \asymp \frac{1}{|x-y|^d\Psi(|x-y|)}$ for all $x,y\in \R^d\, .$
It is now straightforward to check that all the conditions of \cite[Theorem 1.2]{CK08} (as well as the geometric condition of 
\cite{CK08})
are satisfied.  
Let
$$
\wt{q}(t,x,y):=\Psi(t)^{-d}\wedge \frac{t}{|x-y|^d\Psi(|x-y|)}\, , \quad t>0, \ x,y\in \R^d.
$$
It follows from 
\cite{CK08} that there exists  a conservative Feller and strongly Feller process $Z$ associated with 
$(\CC,  \DD(\CC))$ 
that can start from every point in $\R^d$. Moreover, the process $Z$ has a continuous transition density $p(t,x,y)$  on $(0, \infty)  \times \R^d\times \R^d$ (with respect to the Lebesgue measure) which satisfies the following estimates: There exists  
$c_1\ge 1$ such that
$$
c_1^{-1}\wt{q}(t,x,y)\le p(t,x,y)\le c_1 \wt{q}(t,x,y)\, ,\quad t>0, \ x,y\in \R^d\, .
$$

Denote the part of the process $Z$ killed upon exiting $U$ by $Z^U$. 
Then the Dirichlet form of $Z^U$ is $(\CC, \DD_U(\CC))$ where $\DD_U(\CC)=\{u\in \DD(\CC):\, u=0 \textrm{ q.e.~on  } \R^d\setminus U\}$. 
By \cite[Theorem 3.3.9]{CF}, $C_c^{\infty}(U)$ is a core of $(\CC, \DD_U(\CC))$.
By the definition of $J_\gamma$, 
we have that for $u,v\in \DD_U(\CC)$,
\begin{eqnarray*}
\CC(u,v)&=&\frac{1}{2}\int_U \int_U (u(x)-u(y))(v(x)-v(y))J_{\gamma}(x,y)\, dy\, dx +\int_U u(x)v(x)\kappa^{Z}_U(x)\, dx\\
&=&\frac{1}{2}\int_U \int_U (u(x)-u(y))(v(x)-v(y))
J(x,y)
\, dy\, dx +\int_U u(x)v(x)\kappa^{Z}_U(x)\, dx
\end{eqnarray*}
with
\begin{equation}\label{e:kappa-Z-new}
\kappa^{Z}_U(x)=\int_{\R^d\setminus U}J_{\gamma}(x,y)\, dy=\gamma \int_{\R^d\setminus U}j(|x-y|)\, dy=\gamma\kappa_U^X(x)\, ,\quad x\in U\, .
\end{equation}

It follows from 
\eqref{e:kappaU-vs-kappaU-X} that
$c_2\kappa_U(x)\le \gamma^{-1}\kappa_U^Z(x)\le c_3\kappa_U(x)$ for all $x\in U$ 
with positive constants 
$c_2$ and $c_3$ 
independent of $\gamma$. 
Let $\gamma =1/c_3$ and fix it. Then with $c_{4}:=\gamma c_2$ we see that
\begin{equation}\label{e:kappa-U-kappa-Z-comparable}
c_4\kappa_U(x)\le \kappa_U^Z(x)\le \kappa_U(x)\, , \qquad x\in U\, .
\end{equation}
It follows that for $u\in C_c^{\infty}(U)$,
\begin{eqnarray*}
\EE^{D,\kappa}_1(u,u)&=&\EE^{D,\kappa}(u,u)+\int_U u(x)^2 \, dx\\
&=&\frac{1}{2}\int_U \int_U (u(x)-u(y))^2 J(x,y)\, dy\, dx +\int_U u(x)^2 \kappa_U(x)\, dx+\int_U u(x)^2 dx \\
&\asymp & \frac{1}{2}\int_U \int_U (u(x)-u(y))^2 J_{\gamma}(x,y)\, dy\, dx +\int_U u(x)^2 \kappa_U^Z(x)\, dx+\int_U u(x)^2 dx\\
& =&\CC(u,u) +\int_U u(x)^2 dx=\CC_1(u,u)\, .
\end{eqnarray*}
Since $C_c^{\infty}(U)$ is a core of both 
$(\EE^{D, \kappa}), \FF^{D,\kappa}_U)$ 
and $(\CC, \DD_U(\CC))$, we conclude that $\FF^{D,\kappa}_U=\DD_U(\CC)$.

We now define $\wt{\kappa}:U\to \R$ by 
\begin{equation}\label{e:wt-kappa}
\wt{\kappa}(x):=\kappa_U(x)-\kappa_U^Z(x),  \quad x\in U.
\end{equation}
By the choice of $\gamma$ we have that $\wt{\kappa}\ge 0$. 
Note that, by \eqref{e:ww1} there exists $c_5>0$ such that 
$$
\kappa_U^Z(x)=
\gamma \int_{\R^d\setminus U}j(|x-y|)\, dy \le \gamma \int_{\R^d\setminus B(x, \delta_U(x))}j(|x-y|)\, dy
\le c_{5}\frac{1}{\Psi(\delta_U(x))},
\quad x\in U.
$$
Hence it follows from \eqref{e:kappa-U-kappa-Z-comparable} that 
\begin{equation}\label{e:estimate-kappaU}
\kappa_U(x)\le 
c_{4}^{-1}\kappa_U^Z(x)\le \frac{c_{4}^{-1}c_{5}}{\Psi(\delta_U(x))}\, , \quad x\in U.
\end{equation}
Let $\mu(dx)=\wt{\kappa}(x)\, dx$ be a measure on $U$. For $t>0$ and $a\ge 0$, define
$$
N^{U, \mu}_a(t):=\sup_{x\in \R^d}\int_0^t \int_{z\in U: \delta_U(z)>a\Psi^{-1}(t)}\wt{q}(s,x,z)\mu(dz)\, ds\, .
$$
By the definition of $\wt{q}$ and 
\eqref{e:estimate-kappaU} one can  check that $\sup_{t<1}N^{U, \mu}_a(t)<\infty$ and $\lim_{t\to 0} N^{V,\mu}_0(t)=0$ for any relatively compact open 
set $V\subset U$, that is, $\mu\in \textbf{K}_1(U)$ in the notation of \cite[Definition 2.12]{CKSV19}. 

Let $A_t:=\int_0^t \wt{\kappa} (Z^U_s)\, ds$. Then $(A_t)_{t\ge 0}$ is a positive continuous additive functional of $Z^U$ in the strict sense (i.e.~without an exceptional set) with Revuz measure $\wt{\kappa} (x)dx$. For any non-negative Borel function $f$ on $U$, let
$$
T^{U, \wt{\kappa} }_t f(x):=\E_x[\exp(-A_t)f(Z^U_t)]\, ,\qquad t>0, x\in U\, ,
$$
be the Feynman-Kac semigroup of $Z^U$ associated with $\wt{\kappa} (x)dx$. By \cite[
Proposition
 2.14]{CKSV19}, the Hunt process $Z^{U, \wt{\kappa} }$ on $U$ corresponding to the transition semigroup $(T^{U,\wt{\kappa} }_t)_{t\ge 0}$ has a transition density $q^U(t,x,y)$ (with respect to the Lebesgue measure) such that $q^U(t,x,y)\le c_{17} \wt{q}(t,x,y)$ for $t<1$. Further, $(t,y)\mapsto q^U(t,x,y)$ is continuous for each $x\in U$.

According to \cite[Theorem 6.1.2]{FOT}, the Dirichlet form $\CC^{U, \wt{\kappa} }$ corresponding to $T^{U, \wt{\kappa} }_t$ is regular and is given by
\begin{align*}
\CC^{U,\wt{\kappa} }(u,v)
=\frac{1}{2}\int_U\int_U (u(x)-u(y))(v(x)-v(y)) J(x,y) dy  dx+\int_U u(x)v(x)\kappa_U(x)\, dx
\end{align*}
with the domain
$
\DD^{\wt{\kappa}}_U=\DD_U(\CC)\cap L^2(U, \wt{\kappa}(x)dx)\, .
$
Since 
$(\CC^{U,\wt{\kappa}}, \DD^{\wt{\kappa}}_U)$ 
is regular, the set $\DD^{\wt{\kappa}}_U\cap C_c(U)=\DD_U(\CC)\cap C_c(U)$ is its core.
By comparing with \eqref{e:EDkappa} we see that 
$$
\EE^{D, \kappa}(u,v)=\CC^{U, \wt{\kappa}}(u,v)\, , \qquad u,v\in C_c^{\infty}(U)\, .
$$
Now we show 
that the Dirichlet spaces 
 $(\EE^{D, \kappa}, \FF^{D,\kappa}_U)$ and $(\CC^{U, \wt{\kappa}}, \DD^{\wt{\kappa}}_U)$ 
 are equal. We know that $C_c^{\infty}(U)$ is a core for $\EE^{D, \kappa}$. 
One can easily check that this is  also true for $\CC^{U, \wt{\kappa}}$.
 Further,  $C_c^{\infty}(U)\subset  C_c(U)\cap \{u\in L^2(U, dx):\, \CC^U(u,u)<\infty\}$ (which is a core). Clearly, $C_c^{\infty}(U)$ is dense in $C_c(U)$ with uniform norm. 
It is easy to see that  $C_c^{\infty}(U)$ is dense in $C_c(U)\cap \{u\in L^2(U, dx):\, \CC^U(u,u)<\infty\}$ with $\CC^{U,\kappa_U}_1$ norm.
Thus the process $Z^{U, \wt{\kappa}}$ coincides with $Y^{\kappa, U}$. \qed

\begin{prop}\label{p:YDregularized}
The process $Y^{ \kappa}$ can be refined to start from every point in $D$.
\end{prop}
\pf 
Using Lemma \ref{l:YDkappaU}, the proof is the same as that of \cite[Proposition 3.2]{KSV21a}.
\qed


\section{Analysis of the generator}\label{s:generator}

In this section we assume that \textbf{(H1)}-\textbf{(H4)} and \eqref{e:assumption-on-j}-\eqref{e:H-infty} and \eqref{e:bound-on-kappa}  hold. 
Let 
$$
C_c^2(D; \R^d)=\{f:D\to  \R : \text{ there exists } u \in  C_c^2(\R^d) \text{ such that } u=f \text{ on } D\}
$$
be the space of functions on $D$ that are restrictions of $C_c^2(\R^d)$ functions. Clearly, if $f\in C_c^2(D; \R^d)$, then $f\in C^2_b(D) \cap L^2(D)$.

For $\epsilon>0$, let
$$
L_{\epsilon}^\sB f(x):=\int_{D, |y-x|> \epsilon}(f(y)-f(x))J(x,y)\, dy- \kappa (x)f(x).
$$
We introduce the operator
\begin{align}\label{e:defn-LBD}
L^\sB f(x):=\textrm{p.v.}\int_{D}(f(y)-f(x))J(x,y)\, dy- \kappa (x)f(x)=\lim_{\epsilon \downarrow 0}L_{\epsilon}^\sB f(x)\, ,\quad x\in D\, ,
\end{align}
defined for all functions $f:D\to \R$ for which the principal value integral makes sense. 
We will show that this is the case when $f\in C_c^2(D; \R^d)$. We start with the following result.

\begin{lemma}\label{l:operator-interpretation-B}
There exists a constant 
$C_{13}>0$ 
such that for any 
bounded 
Lipschitz function $f$ with Lipschitz constant $L$, 
any $x\in D$ and any $r\in (0,\delta_D(x)]$,
\begin{align}\label{e:Lemma 3.3}
& \int_D \left|f(y)-f(x)\right|j(|y-x|) |\sB(x,x)-\sB(x,y)|\, dy
\le C_{13}\Psi(r)^{-1}\Big(
\|f\|_{\infty} +{rL} \Big)\, .
\end{align}
\end{lemma}
\pf
First note that
\begin{align*}
&\int_D \left|f(y)-f(x)\right|j(|y-x|) |\sB(x,x)-\sB(x,y)|\, dy \\
& \le \int_{D, |y-x|<r/2}\left|f(y)-f(x)\right|j(|y-x|) |\sB(x,x)-\sB(x,y)|\, dy\\
&  +\sB(x,x)\int_{D, |y-x|\ge r/2}\left|f(y)-f(x)\right|j(|y-x|) \, dy + \int_{D, |y-x|\ge r/2}\left|f(y)-f(x)\right|j(|y-x|)\sB(x,y)\, dy\\
& =:I_1+I_2+I_3.
\end{align*}
It follows from $\delta_D(x)\ge r$ that, if $|y-x|<r/2$, then $\delta_D(y)>r/2$ and thus
$\delta_D(y)\wedge \delta_D(x)>r/2>|y-x|$. 
 Hence, when $\delta_2  \ge 1/2$,
 by \textbf{(H4)}, 
 \eqref{e:assumption-on-j} and \eqref{e:H-infty},
\begin{eqnarray*}
I_1 &\le &  
C_5 L\int_{D, |y-x|<r/2} |x-y|j(|x-y|)\left(\frac{|x-y|}{\delta_D(x)\wedge \delta_D(y)}\right)^{\theta}\, dy\\
& \le & 
c_1 L  r^{-\theta} \int_{|y-x|<r/2}|y-x|^{1+\theta}|y-x|^{-d}\Psi(|y-x|)^{-1}\, dy \\
&\le & c_2  
L r^{-\theta}\frac{1}{\Psi(r)}\int_0^{r/2}\frac{s^{\theta}\Psi(r)}{\Psi(s)}\, ds
\\&\le & a_2 c_2 
L r^{-\theta}\frac{r^{2\delta_2}}{\Psi(r)}\int_0^{r/2}s^{\theta-2\delta_2}\, ds \le   
c_3 L r\Psi(r)^{-1}\, .
\end{eqnarray*}
When 
$\delta_2 <1/2$
by \textbf{(H2)}, 
 \eqref{e:assumption-on-j} and \eqref{e:H-infty},
\begin{align*}
&I_1\le
2C_1 L\int_{D, |y-x|<r/2} |x-y|j(|x-y|)\, dy\\
& \le
c_4 L   \int_{|y-x|<r/2}|y-x|^{-d+1}\Psi(|y-x|)^{-1}\, dy 
\le c_5
L \frac{1}{\Psi(r)}\int_0^{r/2}\frac{\Psi(r)}{\Psi(s)}\, ds \\
&\le a_2 c_5
L\frac{r^{2\delta_2}}{\Psi(r)}\int_0^{r/2}s^{-2\delta_2}\, ds \le   
c_6 L r\Psi(r)^{-1}\, .
\end{align*}

Next,
\begin{align*}
& I_2\le 2\|f\|_{\infty}\sB(x,x) \int_{|y-x|>r}j(|x-y|)\, dy\le  c_7\|f\|_{\infty}
 \int_{r/2}^{\infty}t^{d-1}j(t)\, dt  \le c_8\|f\|_{\infty}
\Psi(r)^{-1}\, .
\end{align*}
Finally, by Lemma \ref{l:H3-strong} 
$$
I_3=\int_{D, |y-x|\ge r/2}\left|f(y)-f(x)\right|J(x,y)\, dy\le 2\|f\|_{\infty}\int_{D, |y-x|\ge r/2} J(x,y)\, dy
\le c_9\|f\|_{\infty}\Psi(r)^{-1}
 \, .
$$
Combining the estimates for $I_1$, $I_2$ and $I_3$ we get \eqref{e:Lemma 3.3}. \qed

For notational convenience, we use $ L^\sB_0 f(x)=L^\sB f(x)$ below. 

\begin{prop}\label{p:operator-interpretation}
(a)
If $f\in C_c^2(D; \R^d)$, then
$L^\sB f$ is well defined for all  $x \in D$ and
$r>0$. For $0\le\epsilon \le r \wedge (\delta_D(x)/2)$, it holds that
\begin{eqnarray}\label{e:operator-interpretation}
L^\sB_\epsilon f(x)&=&\sB(x,x) \int_{y \in \R^d, \, |x-y| \ge \epsilon}{(u(y)-u(x)-\nabla u(x) {\bf 1}_{\{ |y-x|<r\}}\cdot (y-x))}j(|x-y|)dy\nn\\
& &+\sB(x,x)\int_{\R^d \setminus D}(u(x)-u(y))j(|x-y|)dy\nn\\
& & +\int_{y\in D, \, |x-y| \ge \epsilon}{(u(y)-u(x))}j(|x-y|)(\sB(y, x)-\sB(x, x))dy
- \kappa(x)u(x),
\end{eqnarray}
where $u\in C_c^2(\R^d)$ is any function such that $u=f$ on $D$.

\noindent (b)
There exists a constant 
$C_{14}>0$ 
such that for any 
$f\in C_c^2(D; \R^d)$, any $x\in D$ and any $r\in (0,\delta_D(x)]$ we have
\begin{align}\label{e:LB-on-u-estimate-b}
&\sup_{0\le\epsilon \le r \wedge (\delta_D(x)/2)}|L_\epsilon^\sB f (x)| \le
C_{14} 
\left(r^2\|\partial^2 u\|_{\infty}+r\|\nabla u\|_{\infty}+\
\|u\|_{\infty}\right)\Psi(r)^{-1}
\end{align}
where $u\in C_c^2(\R^d)$ is any function such that $u=f$ on $D$.

\noindent
(c) 
There exists 
$C_{15}>0$ 
such that for any $f\in C_c^2(D; \R^d)$,
any open $U\subset D$ and any $0<r \le 
\delta_U/2$, 
\begin{align}\label{e:LB-on-u-estimate-d}
 &
 \sup_{0\le\epsilon \le r  }\| (L^\sB_\epsilon f)_{|U} \|_{\infty} 
\le  C_{15} \left(
r^2\|\partial^2 u\|_{\infty}+
r\|\nabla u\|_{\infty}+\|u\|_{\infty}\right) \Psi(
r)^{-1}\, ,
\end{align}
where $u\in C_c^2(\R^d)$ is any function such that $u=f$ on $D$.
\end{prop}

\begin{remark}\label{r:operator-interpretation}{\rm 
We note that the value of the right-hand side of \eqref{e:operator-interpretation} does not depend on the choice of $u\in C_c^2(\R^d)$ such that $u=f$ on $D$. This will be seen from the proof below. On the other hand, 
the quantities $\|\partial^2 u\|_{\infty}$, $\|\nabla u\|_{\infty}$, $\
\|u\|_{\infty}$ on the
right-hand sides in \eqref{e:LB-on-u-estimate-b}--\eqref{e:LB-on-u-estimate-d} depend on the choice of $u$, but this is inconsequential for our purpose.
}
\end{remark}
\noindent
\textbf{Proof of Proposition \ref{p:operator-interpretation}:}  (a)
Using Lemma \ref{l:operator-interpretation-B}, the proof is the same as that of \cite[Proposition 3.4(a)]{KSV21a}.

We give the proof for reader's convenience. 
 Let $u\in C_c^2(\R^d)$ be such that $u=f$ on $D$. Fix $x\in D$ and let 
  $\epsilon<r \wedge (\delta_D(x)/2)$. Then
\begin{align}
&\int_{D, \, |x-y|>\epsilon} {(f(y)-f(x))}j(|x-y|)\sB(x,y) dy \nn\\
&=\sB(x, x)\int_{D, \, |x-y|>\epsilon} {(u(y)-u(x))}j(|x-y|)dy\nn\\
&\qquad +\int_{D, \, |x-y|>\epsilon} {(u(y)-u(x))}j(|x-y|)(\sB(x,y)-\sB(x, x))dy\nn\\
&=\sB(x, x)\int_{|x-y|>\epsilon} {(u(y)-u(x))}j(|x-y|)dy+
\sB(x, x)\int_{\R^d \setminus D, \,  |x-y|>\epsilon}{(u(x)-u(y))}j(|x-y|)dy\nn\\
&\qquad +\int_{D, \, |x-y|>\epsilon}{(u(y)-u(x))}j(|x-y|)(\sB(x,y)-\sB(x, x))dy\nn\\
&=\sB(x, x)\int_{|x-y|>\epsilon}{(u(y)-u(x)-\nabla u(x) {\bf 1}_{\{ |y-x|<r\}}\cdot (y-x))}j(|x-y|)dy\nn\\
&\qquad +\sB(x, x)\int_{\R^d \setminus D, \, |x-y|>\epsilon}(u(x)-u(y))j(|x-y|)dy\nn\\
&\qquad +\int_{D, \, |x-y|>\epsilon}{(u(y)-u(x))}j(|x-y|)(\sB(y, x)-\sB(x, x))dy. \label{e:additional-label}
\end{align}
In the last integral above, we have used \textbf{(H1)}.
By subtracting $\kappa(x)u(x)$, we see that \eqref{e:operator-interpretation} holds true, and that the right-hand side of the equality does not depend on the particular choice of the function $u$.
By letting $\epsilon \to 0$ in \eqref{e:additional-label}  and using 
Lemma \ref{l:operator-interpretation-B} (with $r$ there being $\delta_D(x)$)
for the third integral, we see that $L^{\sB}f$ is well defined.

\noindent 
(b) 
Let $u\in C_c^2(\R^d)$ 
be any function such that $u=f$ on $D$. 
Fix $x\in D$ and let $r\in (0, \delta_D(x)]$ and $0\le\epsilon \le r \wedge (\delta_D(x)/2)$. Then by part (a), 
\begin{align*}
L^\sB_\epsilon f(x)&= \sB(x,x)\int_{y \in \R^d, \, |x-y| \ge \epsilon}\left(u(y)-u(x)-\nabla u(x){\bf 1}_{|y-x|<r}\cdot (y-x)\right)j(|y-x|)\, dy \\
& \quad +\sB(x,x)\int_{\R^d\setminus D} (u(x)-u(y))j(|y-z|)\, dy\\
&\quad  +\int_{y \in D, \, |x-y| \ge \epsilon}(f(y)-f(x))j(|y-x|)(\sB(y,x)-\sB(x,x))\, dy -\kappa(x)f(x)\\
&=:I_\epsilon+II+III_\epsilon+IV\, .
\end{align*}
For $
I_\epsilon
$, we use
$$
\left|u(y)-u(x)-\nabla u(x){\bf 1}_{|y-x|<r}\cdot (y-x)\right|\le \|\partial^2 u\|_{\infty}|y-x|^2{\bf 1}_{|y-x|\le r}+2\|u\|_{\infty}{\bf 1}_{|y-z|\ge r}
$$
to get 
\begin{align*}
\sup_{0\le\epsilon \le r \wedge (\delta_D(x)/2)}|I_\epsilon|
&\le \sB(x,x)\int_{\R^d}\left(\|\partial^2 u\|_{\infty}|y-x|^2 {\bf 1}_{|y-x|\le r}+2\|u\|_{\infty}{\bf 1}_{|y-x|\ge r}\right)j(|y-x|)\, dy\\
&\le 
c_1
\left(\|\partial^2 u\|_{\infty}\int_0^r t^{d-1}t^2 t^{-d}\Psi(t)^{-1}\, dt+\int_r^{\infty}2\|u\|_{\infty}t^{d-1}t^{-d}\Psi(t)^{-1}\, dt\right)\\
&\le 
c_2
(\|\partial^2 u\|_{\infty}r^2+ 2\|u\|_{\infty})\Psi(r)^{-1}\, .
\end{align*}
For $II$ we use  $\delta_D(x)\ge r$ 
to get 
\begin{align*}
|II|
\le 2\sB(x,x)\|u\|_{\infty}\int_{B(x, \delta_D(x))}j(|y-x|)\, dy
 \le 
c_3\|u\|_{\infty} \Psi(\delta_D(x))^{-1}\le 
c_3\|u\|_{\infty} \Psi(r)^{-1}\, .
\end{align*}
$\sup_{0\le\epsilon \le r \wedge (\delta_D(x)/2)}|III_\epsilon|$
is estimated in  
Lemma \ref{l:operator-interpretation-B} (with $L=\|\nabla u\|_{\infty}$), 
while for $IV$ we use \eqref{e:bound-on-kappa} to get
$$
|IV|
\le C_1\|f\|_{\infty}\Psi(\delta_D(x))^{-1}\le C_1\|f\|_{\infty} \Psi(r)^{-1}\, .
$$

\noindent
(c) 
Recall that  $
r \le \delta_D(x)$ for $x\in U$. 
Thus, using 
\eqref{e:LB-on-u-estimate-b},
\begin{align*}
&\sup_{0\le\epsilon \le r  }\| (L^\sB_\epsilon f)_{|U} \|_{\infty}  \le 
c_4 \left(r^2\|\partial^2 u\|_{\infty}+
r\|\nabla u\|_{\infty}+\|u\|_{\infty}\right) 
\Psi(
r)^{-1}.
\end{align*}
\qed

\begin{corollary}\label{c:L2-generator}
 Let $(A, \DD(A))$ be the $L^2$-generator of the semigroup 
corresponding to ${\mathcal E}^{D, \kappa}$.
Then  $C_c^2(D; \R^d)\subset \DD(A)$ and 
$A|_{C_c^2(D; \R^d)}=L^\sB|_{C_c^2(D; \R^d)}$. 
\end{corollary}
\pf Since $\kappa$ is locally bounded, it suffices to show that, for $u\in  C_c^2(D; \R^d)$ and  $v \in C_c^2(D)$, 
\begin{align}
\label{e:DFY_1}
& \int_{D}\int_{D}{(u(y)-u(x))(v(y)-v(x))}J(x,y)(x,y) \, dy\, dx=
2\int_{D}\big(L^\sB u(x)-\kappa(x)u(x)\big) v(x)\, dx.
\end{align}
By Lemma \ref{l:df-well-defined}, the left-hand side is well defined and absolutely integrable.
Hence by the dominated convergence theorem and the symmetry of $\sB$, 
\begin{align}
\label{e:DFY_2}
& \int_{D}\int_{D} {(u(y)-u(x))(v(y)-v(x))}j(|x-y|)\sB(x,y) dydx
\nn\\
&=
\lim_{\epsilon \downarrow
0}  
\int_{D}
\int_{y \in D:|x-y|>\epsilon} {(u(y)-u(x))(v(y)-v(x))}j(|x-y|)\sB(x,y) dydx\nn\\
&=
2\lim_{\epsilon \downarrow
0} 
 \int_{D}
\int_{y \in D:|x-y|>\epsilon} {(u(y)-u(x))}j(|x-y|)\sB(x,y) dy\,  v(x)dx \nn\\
&=2\lim_{\epsilon \downarrow 0} \int_{\text{supp}(v)}\int_{y \in D:|x-y|>\epsilon} {(u(y)-u(x))}j(|x-y|)\sB(x,y) dy\,  v(x)dx.
\end{align}
Let $\epsilon < \epsilon_0:=\mathrm{dist}(\partial D,\mathrm{supp}(v))/2$. 
It follows from Proposition \ref{p:operator-interpretation} (c) 
(by taking $U=\mathrm{supp}(v)$) that  
\begin{align*}
&\sup_{x \in \mathrm{supp}(v),\epsilon < \epsilon_0}
\left|\int_{y \in D:|x-y|>\epsilon} {(u(y)-u(x))}j(|x-y|)\sB(x,y) dy  \right| \\
&\qquad \le c_1 \left(\epsilon_0^2\|\partial^2 u\|_{\infty}
+\epsilon_0\|\nabla u\|_{\infty}+\|u\|_{\infty}\right) \Psi(\epsilon_0)^{-1}.
\end{align*}
Since the right-hand side is finite,  we can use the dominated convergence theorem to conclude that \eqref{e:DFY_1} holds. \qed

Corollary \ref{c:L2-generator} says that $L^\sB$ is the extended generator of the semigroup $(T_t)_{t\ge0}$ corresponding to ${\mathcal E}^{D, \kappa}$.

\medskip  
Let $U$ be an open set with $U\subset \overline{U}\subset D$. 
Recall that $\kappa^U$  and $\kappa_U=\kappa+\kappa^U$ are defined in 
\eqref{e:kappaU}.
Consider now the process $Y^{\kappa, U}$.
Denote $L^\sB_U u:=L^{\sB,U} u-\kappa^U(\cdot) u$, where
$$
L^{\sB,U} u(z):=\mathrm{p.v.} \int_U (u(y)-u(z))J(y,z)\, dy -\kappa(z)u(z)\, ,\qquad u\in U.
$$
Since $\kappa_U=\kappa+\kappa^U$, we can write
$$
L^{\sB}_{U} u(z)=\mathrm{p.v.} \int_U (u(y)-u(z))J(y,z)\, dy -\kappa_U(z)u(z)\, ,\qquad u\in U\, .
$$

\begin{corollary}\label{c:L2-generatork}
Let $U$ be an open subset of $D$ and let $(A, 
\DD(A))$ be the $L^2$-generator of the semigroup $(T_t)_{t\ge0}$ of $Y^{\kappa, U}$.
Then  $C_c^2(U)\subset \DD(A)$ and 
$A_{|C_c^2(U)}=(L^\sB_U)_{|C_c^2(U)}=L^\sB_{|C_c^2(U)}$.
\end{corollary}
\pf 
If $u\in C_c^2(U)$, then for $z\in U$, 
\begin{align*}
&L^\sB u(z)=\lim_{\epsilon \to 0}\int_{D, |y-z|<\epsilon}(u(y)-u(z))J(y,z)dy -\kappa(z)u(z)\\
&=\lim_{\epsilon \to 0}\int_{U, |y-z|<\epsilon}(u(y)-u(z))J(y,z)dy+\lim_{\epsilon \to 0}\int_{D\setminus U, |y-z|<\epsilon}(u(y)-u(z))J(y,z)dy-\kappa(z)u(z)\\
&=L^{\sB,U} u(z)-\kappa^U(z)u(z)-\kappa(z)u(z)= L^\sB_U u(z)\, .
\end{align*}
Thus, the corollary follows from Corollary \ref{c:L2-generator} and its proof.
\qed

The goal of the remainder of this section is 
to  prove
a Dynkin-type formula (Theorem  \ref{t:Dynkin}), which will be used in  \cite{KSV22b}.

Recall that, for an open set $U\subset D$, $\tau_U=\tau_U^{Y^{\kappa}}=\inf\{t>0:\, Y^{\kappa}_t\notin U\}$.
\begin{lemma}\label{l:martingale-U}
Suppose that $U$ is an open set with $U\subset \overline{U}\subset D$.
For any $u\in C_c^2(D)$ and any $x\in U$,
\begin{equation}\label{e:martingale-U}
M_t^u:=u(Y^{\kappa, U}_t)-u(Y^{\kappa, U}_0)-\int_0^t 
L^\sB u(Y^{\kappa, U}_s)\, ds
\end{equation}
is a $\P_x$-martingale  with respect to the filtration of $Y^{\kappa, U}$.
\end{lemma}
\pf 
We first assume that $u\in C_c^2(U)$ and we follow the proof of \cite[Lemma 2.2]{GKK}. 
Let $(A, \DD(A))$ be the $L^2$-generator of the semigroup $(T_t)_{t\ge 0}$ of $Y^{\kappa, U}$. By Corollary \ref{c:L2-generatork},
$C_c^2(U)\subset \DD(A)$ and 
$A_{|C_c^2(U)}=(L^\sB_U)_{|C_c^2(U)}=L^\sB_{|C_c^2(U)}$. 
Since
$\|(T_tf-f)-\int^t_0T_sAfds\|_{L^2(U)}=0$, 
\begin{align}
\label{DF1}
T_t u(x)-u(x)=\int_0^t T_s 
L^\sB
 u(x)\, ds\,  \qquad \textrm{a.e. } x\in U\, . 
\end{align}
Let $g_t(x):=\int_0^t T_s 
L^\sB
u(x)\, ds$, $x\in U$. 
Note that  $L^\sB u$ is bounded in $U$ by Proposition \ref{p:operator-interpretation} (c).
Thus,  $|g_t(x)|\leq t
||(L^\sB u)_{|U}||_\infty<\infty$ for all $x\in U$. 
Since  $Y^{\kappa, U}$ is strongly Feller by Lemma \ref{l:YDkappaU},,
we have $T_\epsilon g_{t-\epsilon}\in C_b(U)$ for all $\epsilon\in (0, t)$. 
Moreover,
$$|g_t(x)-T_\epsilon g_{t-\epsilon}(x)|=|g_\epsilon(x)|\leq \epsilon 
||(L^\sB u)_{|U}||_\infty, \quad \text{for all }x\in U.$$
Hence, $g_t$ is continuous and  \eqref{DF1} holds for any $x\in U$. 
Using this and 
the Markov property,
we get the desired conclusion for $u\in C_c^2(U)$.

In general, when $u\in C_c^2(D)$, let $V$ be a compact open subset of $D$ such that supp$(u) \cup \overline U \subset V$.
By the conclusion above, we have that for any $u\in C_c^2(D)$ and any $x\in U$,
$$
u(Y^{\kappa, V}_t)-u(Y^{\kappa, V}_0)-\int_0^t 
L^\sB u(Y^{\kappa, V}_s)\, ds=u(Y^{\kappa}_t)\mathbf{1}_{t<\tau_V}-u(Y^{\kappa}_0)-\int_0^{t\wedge \tau_V}L^\sB u(Y^{\kappa}_s)ds
$$
is a $\P_x$-martingale  with respect to the filtration of $Y^{\kappa}$. Since 
$\tau_U \le \tau_V$, by the optional stopping theorem we get the desired conclusion for $u\in C_c^2(D)$.
\qed

\begin{prop}\label{p:martingale-U} 
Suppose that $U$ is an open set with $U\subset \overline{U}\subset D$.
For any $u\in C_c^2(D)$ and any $x\in U$,
\begin{equation}\label{e:martingale_U2}
M_t^u=u(Y^{\kappa}_t)\mathbf{1}_{t<\tau_U}-u(Y^{\kappa}_0)-\int_0^{t\wedge \tau_U}L^\sB u(Y^{\kappa}_s)ds\, 
\end{equation}
is a $\P_x$-martingale  with respect to the filtration of $Y^{\kappa}$.
\end{prop}
\pf 
Note that 
$L^\sB u(Y^{\kappa}_s)\mathbf{1}_{s<\tau_U\wedge \zeta}
=L^\sB u(Y^{\kappa}_s)\mathbf{1}_{s<\tau_U}$ and that $u(Y^{\kappa, U}_t)=u(Y^{\kappa}_t)\mathbf{1}_{t<\tau_U}$. 
Thus we can rewrite \eqref{e:martingale-U} as \eqref{e:martingale_U2}. \qed

For any $x\in D$ and Borel subset $A$ of $D_{\partial}$, we define
$N(x, A)=\int_{A\cap D}J(x, y)dy+\kappa(x){\bf 1}_A(\partial)$.  
Then it is known that $(N, t)$ is a L\'evy system for $Y^{\kappa}$ (cf.~\cite[Theorem 5.3.1]{FOT} and the argument in \cite[p.40]{CK03}), that is, for any non-negative Borel function $f$ on $D\times D_{\partial}$ vanishing on the diagonal and any stopping time $T$,
\begin{align}\label{e:levys}
\E_x\sum_{s\le T}f(Y^{\kappa}_{s-}, Y^{\kappa}_s)=
\E_x\left(\int^T_0\int_{D_{\partial}}f(Y^{\kappa}_s, y)N(Y^{\kappa}_s, dy)ds\right), \quad x\in D.
\end{align}

We are now ready to establish
 the following Dynkin-type theorem.
 
 \begin{thm}\label{t:Dynkin}
 Suppose that $U$ is an open set with $U\subset \overline{U}\subset D$.
For any  non-negative function $u$ defined on $D$ satisfying 
$u\in C^2(\overline{U})$ and any $x\in U$,
\begin{equation}\label{e:Dynkin}
\E_x[u(Y^{\kappa}_{\tau_U})]=u(x)+\E_x\int_0^{\tau_U}L^\sB u(Y^{\kappa}_s)ds\, .
\end{equation}
 \end{thm}
\pf
For any  non-negative function $u$ on $D$ satisfying
$u\in C^2(\overline{U})$, choose an open set $V$ of $D$ and
$f\in C_c^2(D)$ such that $U \subset \overline{U} \subset V \subset \overline{V} \subset D$,
and $f=u$ on $V$  and $f\le u$ on $D$. 
Let  $h:=u-f$  so that $u=h+f$,  $h\ge 0$,   and $h=0$ on $V$.
Since $f\in C_c^2(D)$, by Proposition \ref{p:martingale-U},  
$$
\E_x[f(Y^{\kappa}_t)\mathbf{1}_{t<\tau_U}]=f(x)+\E_x\int_0^{t\wedge \tau_U}L^\sB f(Y^{\kappa}_s)ds\, .
$$
Proposition \ref{p:operator-interpretation}(c) implies that 
$||(L^\sB f)_U ||_\infty<\infty$. 
Thus, by letting $t \to \infty$, 
\begin{equation}\label{e:Dynkin1}
\E_x[f(Y^{\kappa}_{\tau_U})]=f(x)+\E_x\int_0^{\tau_U}L^\sB f(Y^{\kappa}_s)ds\, .
\end{equation}
On the other hand, 
since $h=0$ on $V$, for $y \in U$,
$$
L^\sB h(y)=
\textrm{p.v.}\int_{D}(h(z)-h(y))J(y,z)\, dz- \kappa (y)h(y)
=
 \int_{D\setminus V} h(z) J(y, z)dz.
$$
Thus, by the L\'evy system formula \eqref{e:levys}
\begin{align}\label{e:Dynkin2}
&\E_x[h(Y^{\kappa}_{\tau_U})]=
\E_x[h(Y^{\kappa}_{\tau_U}): Y^{\kappa}_{\tau_U} \in 
D \setminus V ]\nn\\
&=
\E_x \int_0^{\tau_U}
 \int_{D\setminus V} h(z) J(Y^{\kappa}_s, z)dz=
\E_x\int_0^{\tau_U}L^\sB h(Y^{\kappa}_s)ds\, .
\end{align}
Adding  \eqref{e:Dynkin1} and \eqref{e:Dynkin2}, we get \eqref{e:Dynkin}.
\qed


\section{Harnack inequality and H\"older continuity of Harmonic functions}\label{s:hi}

In this section we assume that 
\textbf{(H1)}-\textbf{(H4)},
\eqref{e:assumption-on-j}-\eqref{e:H-infty} and \eqref{e:bound-on-kappa} hold.

\begin{lemma}\label{l: exit-time-probability}
There exists a constant 
$C_{16}>0$ 
such that for all 
$x\in D$ and $r>0$ with  
$B(x, 2r )\subset D$, 
$$
\P_x(\tau_{B(x,r)} < t \wedge \zeta)\le 
C_{16}t\ \Psi(r)^{-1} \, .
$$ 
\end{lemma}

\pf 
Let $x\in D$ and $r>0$ be such that  
$B(x, 2r )\subset D$.
Let $f:\R^d \to [-1, 0]$ be a $C^2$ function such that $f(z)=-1$ for $|z|\le 1/2$, $f(y)=0$ for $|z|\ge 1$ and that $\|\nabla f\|_{\infty}+\|\partial^2 f\|_{\infty}=:c_1<\infty$. Define
$$
f_r(y):=f\left(\frac{y-x}{r}\right)\, .
$$
Then
$f_r\in  C_c^2(D)$, 
$f_r(y)=-1$ for $y\in B(x,r/2)$ and $f_r(y)=0$ for $y\in D\setminus B(x,r)$. 
By \eqref{e:martingale_U2}, 
$$
f_r(Y^{\kappa}_t)\1_{t<\tau_{B(x,r)}}-f_r(Y^{\kappa}_0)-\int_0^{t\wedge \tau_{B(x,r)}} L^\sB f_r (Y^{\kappa}_s)\, ds
$$
is a $\P_y$-martingale for every $y\in B(x,2r)$. 
Hence,
\begin{align}
&
\P_x\big(\tau_{B(x,r)}<t\wedge \zeta\big)=\P
_x\big(|Y^{\kappa}_{\tau_{B(x,r)}\wedge t}-x|\ge r, \tau_{B(x,r)}\wedge t<\zeta\big)\nonumber\\
& =\E
_x\big[1+f_r(Y^{\kappa}_{\tau_{B(x,r)}\wedge t}), |Y^{\kappa}_{\tau_{B(x,r)}\wedge t}-x|\ge r, \tau_{B(x,r)}\wedge t<\zeta\big] \nonumber \\
&\le \E
_x \big[1+f_r(Y^{\kappa}_{\tau_{B(x,r)}\wedge t})\big]=
-f_r(x)
+\E
_x\big[f_r(Y^{\kappa}_{\tau_{B(x,r)}\wedge t})\big]=\E
_x\left[\int_0^{\tau_{B(x,r)}\wedge t} L^\sB f_r(Y^{\kappa}_s)\, ds\right] \nonumber\\
&\le \|(L^\sB f_r)_{|B(x,r)}\|_{\infty}\, \E_x[\tau_{B(x,r)}\wedge t]  \le t \|(L^\sB f_r)_{|B(x,r)}\|_{\infty}
  \, .\label{e:PE}
\end{align}
The first inequality 
follows because $1+f_r\ge 0$. Note that here $f_r(Y^{\kappa}_{\tau_{B(x,r)}\wedge t})$ makes sense regardless whether $\tau_{B(x,r)}\wedge t<\zeta$ or not (by definition $f_r(\partial)=0$). 
Since $\| f_r\|_{\infty}+r\|\nabla f_r\|_{\infty}+r^2\|\partial^2 f_r\|_{\infty}=1+c_1$, applying Proposition \ref{p:operator-interpretation} (c),
we get the desired conclusion.
\qed

\begin{lemma}\label{l:new}
For all $x\in D$ and all $r>0$ with  $B(x,2r)\subset D$, it holds that
$\P_x(\tau_{B(x,r)}=\zeta < t)\le C_{4}\Psi(r)^{-1}t$.
\end{lemma}

\pf
By the L\'evy system formula,
\begin{align*}
\P_x(\tau_{B(x,r)}=\zeta < t)=\E_x\sum_{s<t}{\bf 1}_{B(x, r)\times\{\partial\}}
(Y^{\kappa}_{s-}, Y^{\kappa}_s)
=\E_x\int^t_0{\bf 1}_{B(x, r)}(Y^{\kappa}_s)\kappa(Y^{\kappa}_s)ds.
\end{align*}
Since $\kappa(y) \le C_4/\Psi(\delta_D(y))\le C_4/\Psi(r)$ for $y\in B(x, r)$ 
by \eqref{e:bound-on-kappa}, we immediately get
$
\P_x(\tau_{B(x,r)}=\zeta < t)\le  C_{4}\Psi(r)^{-1}t.
$
\qed

Let $A(x, r_1, r_2)$ denote the annulus $\{y\in \R^d: r_1\le |y-x|<r_2\}$.

\begin{prop}\label{p:exit-time-estimate}
(a) 
 There exists a constant 
$C_{17}>0$ 
such that for all $x_0\in D$ and   $r>0$
 with $B(x_0,r) \subset D$, it holds that
$$
\E_x \tau_{B(x_0,r)}\ge C_{17} \Psi(r)\, ,
\quad x\in B(x_0,r/2).
$$

\noindent
(b) 
For every $\epsilon >0$,  there exists  
$C_{18} =C_{18}(\epsilon) >0$ 
such that for all $x_0\in  D$ and 
$r >0$ satisfying 
$B(x_0, (1+\epsilon)r)\subset D$, it holds that
$$
\E_x \tau_{B(x_0,r)}\le 
C_{18}\Psi(r)\, ,\quad x\in B(x_0,r)\, .
$$
\end{prop}
\pf 
(a) 
Let $x\in D$ and $r>0$ be such that $B(x,r)\subset D$. 
It follows from Lemmas \ref{l: exit-time-probability}--\ref{l:new} and \eqref{e:H-infty}
that
$$
\P_x(\tau_{B(x, r/2)}< t) \le c_{1} \Psi(r)^{-1}t\, .
$$
Therefore,
$$
\E_x \tau_{B(x, r/2)}\ge t \P_x( \tau_{B(x, r/2)}
\ge t)\ge t(1-c_{1}\Psi(r)^{-1}t), \quad t>0.
$$
Choose $t=\Psi(r)/(2c_{1})$, so that $1-c_{1}\Psi(r)^{-1}t=1/2$. Then
$$
\E_x \tau_{B(x, r/2)} \ge \frac12 \Psi(r)/(2c_{1})=c_{2}\Psi(r)\, .
$$
Now let $B(x_0,r)\subset D$ and  
$x\in B(x_0,r/2)$. Then 
$B(x,r/2)\subset B(x_0,r)\subset D$. 
By what was proven above,
$$
\E_x \tau_{B(x_0,r)}\ge \E_x \tau_{B(x,r/2)}\ge c_{2}\Psi(r)\, .
$$

\noindent
(b)
Let $\epsilon_0 :=\epsilon/3$, $x_0\in  D$ and $r>0$ be such that $B(x_0, (1+3\epsilon_0)r)\subset D$.
For $y\in B(x_0,r)$ and $u\in A(x_0, (1+\epsilon_0)r, (1+2\epsilon_0)r)\big)$, 
$\delta_D(u)\wedge \delta_D(y)\ge \epsilon_0 r\ge (\epsilon_0
/(2+2\epsilon_0))|u-y|$. 
Thus, by \textbf{(H2)}, and then using \eqref{e:assumption-on-j}-\eqref{e:H-infty},
 $$
 J(u,y)\ge c_3 j(|u-y|)\ge c_4 j(|u-x_0|), 
 \quad (y,u)\in B(x_0,r) \times A(x_0, (1+\epsilon_0)r, (1+2\epsilon_0)r)\big).$$
Therefore, for $y\in B(x_0,r)$
\begin{align}\label{e:ddsw}
\int_{A(x_0, (1+\epsilon_0)r, (1+2\epsilon_0)r)}J(u,y)\, du\ \ge&
c_4\int_{A(x_0, (1+\epsilon_0)r, (1+2\epsilon_0)r)}j(|u-x_0|)\, du\nn\\
\ge & 
c_5 \int_{(1+\epsilon_0)r}^{(1+2\epsilon_0)r}\frac{1}{t\Psi(t)}\, dt \ge c_6 \frac{1}{\Psi(r)}\, .
\end{align}
 For $x\in B(x_0,r)$, by using \eqref{e:ddsw} in the 
 last inequality below,
\begin{align*}
1\ge &\P_x\big(Y^{\kappa}_{\tau_{B(x_0,r)}} \in A(x_0, (1+\epsilon_0)r, (1+2\epsilon_0)r)\big)\nonumber\\
= &\E_x\int_0^{\tau_{B(x_0,r)}}\int_{A(x_0, (1+\epsilon_0)r, (1+2\epsilon_0)r)}J(u,Y^{\kappa}_s)\, du\, ds
\ge c_6\E_x\tau_{B(x_0,r)}/{\Psi(r)}, 
\end{align*}
which is the required inequality. \qed

Let $T_A$  be the first hitting time to $A$ for $Y^\kappa$. 

\begin{lemma}\label{l:A2}
For every $\epsilon\in (0,1)$ there exists $C_{19}=C_{19}(\epsilon)>0$ such that for 
all $x\in D$ and $r>0$ with $B(x,(1+3\epsilon)r)\subset  D$,
 and any Borel set $A\subset B(x,r)$,
$$
\P_y(T_A<\tau_{B(x,(1+2\epsilon)r)})\ge C_{19}\frac{|A|}{|B(x,r)|}\, , \qquad y\in B(x,(1+\epsilon)r).
$$
\end{lemma}
\pf Without loss of generality we assume that 
${\mathbb P}_y(T_A < {\tau}_{B(x,(1+2\epsilon)r)})<1/4$. 
Set $\tau={\tau}_{B(x,(1+2\epsilon)r)}$. For $y\in B(x,(1+\epsilon)r)$ we have that $B(y,2\epsilon r)\subset  D$ 
and $B(y, \epsilon r)\subset B(x, (1+2\epsilon)r)$.
Hence by Lemmas \ref{l: exit-time-probability} and \ref{l:new}, for any 
$y\in B(x,(1+\epsilon)r)$,
$$
\P_y(\tau< t)\leq \P_y(\tau_{B(y,\epsilon r)}< t) \leq c_0\Psi(\epsilon r)^{-1} t\le c_0 a_2 \epsilon^{-2\delta_2}\Psi(r)^{-1} t
=:c_1 \Psi(r)^{-1} t\, .
$$
 Choose $t_0= \Psi(r)/(4c_1)$, so that $\P_y(\tau< t_0) \leq 1/4$. Further,  
if $z\in B(x, (1+2\epsilon)r)$ and $u\in A \subset B(x,r)$, then $|u-z| \leq 2(1+\epsilon)r$. 
By \eqref{e:assumption-on-j} and  \eqref{e:H-infty}, 
$j(|u-z|) \geq 
 c_2r^{-d}/\Psi(r)$
 for some $c_2=c_2(\epsilon)>0$.
  Moreover, $\delta_D(u)\wedge \delta_D(z)\ge \epsilon r \ge \frac{\epsilon}{2(1+\epsilon)}|u-z|$, 
 implying by \textbf{(H2)} that $\sB(u,z)\ge c_3$  
($c_3=C_1(\epsilon/(2(1+\epsilon))$). 
Thus,
\begin{align*}
&\P_y (T_A < \tau) \geq  \E_y \sum_{s\leq T_A \wedge \tau \wedge t_0}
{\mathbf 1}_{\{Y^{\kappa}_{s-}\neq 
Y^{\kappa}_{s},
Y^{\kappa}_s\in A\}} 
\\
 &=  \E_y \int_0^{T_A \wedge \tau \wedge t_0}
\int_A j(|u-Y^{\kappa}_s|)\sB(u, Y^{\kappa}_s)\, du \, ds 
    \ge   
    \frac{c_2c_3}{r^d \Psi(r)}
    |A| \E_y[T_A \wedge \tau \wedge t_0],
\end{align*}
where in the second line we used properties of the L\'evy system.
Next,
\begin{align*}
{\mathbb E}_y[T_A\wedge \tau \wedge t_0]   
 \ge  t_0 {\mathbb P}_y(T_A \ge \tau \ge t_0) 
 \ge  t_0[1-{\mathbb P}_y(T_A < \tau)-{\mathbb P}_y(\tau <t_0)] 
 \ge  \frac{t_0}{2} = \frac{\Psi(r)}{8 c_1}\, .
\end{align*}
The last two displays give that
$$
{\mathbb P}_y (T_A < \tau) 
 \ge     \frac{c_2c_3}{8 c_1 r^d}|A|= 
 c_4\frac{|A|}{|B(x,r)|}\, , \qquad y\in B(x,(1+\epsilon)r).
$$
 \qed

\begin{lemma}\label{l:A3}
Suppose further that \textbf{(H5)} holds.
There exist 
$C_{20}>0$ and $C_{21}>0$ 
 with the property that if 
$r>0$, and $x\in D$ are such that $B(x,2r)\subset D$,  and $H$ is a bounded
non-negative function with support in $D\setminus B(x,2r)$, then for every $z\in B(x,r)$,
$$
C_{20} \E_z [{\tau}_{B(x,r)}] \int_D H(y)
J(x,y) \, dy \le \E_z H(Y^{\kappa}_{{\tau}_{B(x,r)}}) \le  
C_{21} {\mathbb E}_z [{\tau}_{B(x,r)}] 
\int_D H(y)
J(x,y) \, dy \, .
$$
\end{lemma}

\pf Let $y\in B(x,r)$ and $u\in D\setminus B(x,2r)$. By \eqref{e:H5-for-J}, $J(u,y)\asymp J(x,y)$. Thus using the L\'evy system  we get
\begin{align*}
\E_z \left[ H(Y^{\kappa}_{\tau_{B(x,r)}}) \right]  &= \E_z \int_0^{\tau_{B(x,r)}} \int_{D\setminus B(x,2r)} H(u) J(u, Y^{\kappa}_s)\, du\, ds \\
&\asymp \E_z \int_0^{\tau_{B(x,r)}} \int_{D\setminus B(x,2r)} H(u)  J(u, x)\, du\, ds.
\end{align*}
 \qed

\medskip

\medskip
\noindent
\textbf{Proof of Theorem \ref{t:uhp*}:} 
\noindent
(a) 
Using Proposition \ref{p:exit-time-estimate} and Lemmas \ref{l:A2} and \ref{l:A3} (instead of (A1)--(A3) in \cite{SV04}), the 
proof of (a) is very similar to the proofs of \cite[Theorem 2.2, Theorem 2.4]{SV04}.
We omit the details. 

\noindent
(b) 
By (a) we can and will assume that $L>2$ and $2r < |x_1-x_2|<Lr$. 
For simplicity, let $B_i=B(x_i,r)$, $i=1,2$. 
Then by using harmonicity in the first inequality,  part (a) in the second inequality, and the L\'evy system formula in the second line, we have
\begin{eqnarray}\label{e:hp2-1}
f(x_1)&\ge &\E_{x_1}\big[f(Y^{\kappa}_{\tau_{B_1}}); Y^{\kappa}_{\tau_{B_1}}\in B(x_2, r/2)\big] \ge C_8^{-1} f(x_2) \P_{x_1}\big(Y^{\kappa}_{\tau_{B_1}}\in B(x_2, r/2)\big)\nonumber \\
&=& C_8^{-1} f(x_2)\E_{x_1} \int_0^{\tau_{B_1}}  \int_{B(x_2, r/2)}
J(Y^{\kappa}_s,z)\, dz\, ds\, .
\end{eqnarray}
For $y\in B_1$ and $z\in B(x_2, r/2)$ we have by \eqref{e:H5-for-J} that $J(y,z)\ge c_1J(x_1, z)$. Further, $\delta_D(x_1)\wedge \delta_D(z)\ge r/2 \ge 
(2L+2)^{-1}
|x_1-z|$, hence by \textbf{(H2)}, 
$J(x_1,z)\ge c_2 j(|x_1-z|)$ where $c_2=C_1 (\tfrac{1}{2(L+1)})$.
By inserting this in \eqref{e:hp2-1}, 
and by using Proposition \ref{p:exit-time-estimate} (a), we obtain
\begin{align*}
f(x_1)&\ge 
c_1c_2C_8^{-1} f(x_2) \E_x\tau_{B_1}\int_{B(x_2,r/2)}j(|x_1-z|)\, dz\\ 
&\ge c_3 c_2f(x_2) \Psi(r) \frac{1}{((L+1)r)^d \Psi((L+1)r)}\, |B(x_2, r/2)|\\
 & \ge c_4c_2 f(x_2) L^{-d} \frac{\Psi(r)}{\Psi((L+1)r)}\ge  c_5c_2 f(x_2) L^{-d} L^{-2\delta_2}\, .
\end{align*}
The last inequality follows from \eqref{e:H-infty}. \qed

We now show that 
a non scale invariant Harnack inequality 
holds under
much weaker assumptions than \textbf{(H2)}-\textbf{(H5)}, 
and introduce 
weaker versions of hypotheses \textbf{(H4)}-\textbf{(H5)}:

\medskip
\noindent
\textbf{(H4-w)}  
If $\delta_2 \ge 1/2$, then there exists $\theta>2\delta_2-1$ such that for  any relatively compact open set $U\subset D$ there
exists $C_{22}=C_{22}(U)$ such that 
$$
|\sB(x, x)-\sB(x,y)|\le 
C_{22}|x-y|^{\theta}\,
\quad \text{ for all }x,y\in U.
$$

\medskip
\noindent
\textbf{(H5-w)}  
 For any relatively compact open set $U\subset D$ and any open set $V$ such that $\overline{U}\subset V\subset D$,
there exists $C_{23}=C_{23}(U, V)\ge 1$ 
$$
C_{23}^{-1}\sB(x_1,z)\le \sB(x_2,z) \le C_{23} \sB(x_1, z)\, ,\quad \text{for all }x_1,x_2\in U, \, \,z\in D\setminus V\, .
$$

It is clear that \textbf{(H4)}, respectively \textbf{(H5)}, imply \textbf{(H4-w)}, respectively \textbf{(H5-w)}. 
\medskip

\begin{prop}
[non scale invariant Harnack inequality]\label{t:uhp}
Suppose $D$ is a proper open subset of $\R^d$
and assume that \textbf{(H1)}, \textbf{(H2-w)}-\textbf{(H5-w)}, \eqref{e:assumption-on-j}-\eqref{e:H-infty} and \eqref{e:bound-on-kappa0} hold. 
For any compact set $K$ and open set $U$ with $K\subset U\subset\overline U\subset D$, there exists a constant 
$C_{24}=C_{24}(K, U)>0$ 
such that for any non-negative function $f$ in $D$ which is 
harmonic in $U$ with respect to $Y^\kappa$, we have
$$
f(x)\le C_{24} f(y), \qquad \text{ for all } x, y\in K.
$$
\end{prop}

\pf Using \textbf{(H2-w)}-\textbf{(H5-w)} instead of \textbf{(H2)}-\textbf{(H5)}, non scale invariant versions  (with constants depending on $r$) of 
Propositions \ref{p:operator-interpretation}(c) and \ref{p:exit-time-estimate} and
Lemmas \ref{l:A2} and \ref{l:A3} 
can be proved. 
Proposition \ref{p:exit-time-estimate},  
Lemmas \ref{l:A2} and \ref{l:A3} 
imply that
conditions (A1), (A2) and (A3) of \cite{SV04} are satisfied for the  process $Y^{\kappa}$ with constants depending on $r$.
 Thus we can repeat the proofs of 
\cite[Theorems 2.2 and 2.4]{SV04} to finish 
the proof. 
Note that conservativeness does not play any role. 
We omit the details. \qed

The above Harnack inequality is not scale invariant since the constant $C_{25}$ in the result depends on each $K$ and $U$ there.
The scale invariant version of Harnack inequality is not possible under  \textbf{(H3-w)} since the value of 
integral \eqref{e:H3} depends on the sets $U$ and $V$ there.

\medskip

 By following the arguments of \cite[Theorem 4.9]{SV04}
and \cite[Theorem 4.1]{BL}, we can prove Theorem \ref{t:holder}. 
Note that $r$ is missing  in \cite[Theorem 3.14]{KSV21a} and  \cite[Theorem 4.9]{SV04}.
We give a full proof here for reader's convenience.
Note that 
\textbf{(H5)} is not assumed in Theorem \ref{t:holder}.

\noindent
\textbf{Proof of Theorem \ref{t:holder}:} 
By Lemma \ref{l:A2}, there exists 
$c_1>0$ such that for all 
$(s, x)\in (0, \infty)\times D$ with $B(x, 5s)\subset D$, and any $A\subset B(x, s)$ with $|A|/|B(x, s)|\ge 1/3$,
\begin{equation}\label{e:4.1}
\P_y(T_A<\tau_{B(x, 3s)})\ge c_1, \qquad y\in B(x, 2s).
\end{equation} 
For $y \in B(x,s)$ and $s'>2s$, we have $B(y, s'/2) \subset B(x, s')$.
Thus, using Lemma 
\ref{l:H3-strong}, we have that,   
for $y \in B(x,s)$ and $s'>2s$ with $B(x, 2s')\subset D$,
$$
\int_{D\setminus B(x, s')}J(y, z)dz \le \int_{D\setminus B(y, s'/2)}J(y, z)dz \le C_{11}
\Psi(s'/2)^{-1}.
$$
Using this and Proposition \ref{p:exit-time-estimate} (b), we obtain  
that for $s'>2s$ with $B(x, 2s')\subset D$,
\begin{align*}
\P_y(Y^\kappa_{\tau_{B(x, s)}}\in D\setminus B(x, s'))&=
\E_y \int_0^{\tau_{B(x, s)}} \int_{D\setminus B(x, s')}J(Y^\kappa_t, z)dzdt\nn \\
&\le c_2\Psi(s)/ 
\Psi(s'/2), \qquad y\in B(x, s).
\end{align*}
Thus, by \eqref{e:H-infty},
\begin{equation}\label{e:4.2}
\P_y(Y^\kappa_{\tau_{B(x, s)}}\in D \setminus B(x, s'))\le 
c_3\frac{s^{2\delta_1}}{(s')^{2\delta_1}}, \qquad y\in B(x, s), 
s'>2s.
\end{equation}
Let
$$
\gamma=1-\frac{c_1}4, \qquad \rho=\frac13\wedge\left(\frac{\gamma}2\right)^{1/(2\delta_1)}
\wedge \left(\frac{c_1\gamma^2}{8c_3}\right)^{1/(2\delta_1)}.
$$
Let $x\in B(x_0, r/2)$.  We will show that
\begin{equation}\label{e:4.3}
\sup_{B(x, \rho^kr)}f-\inf_{B(x, \rho^kr)}f\le 
\|f\|_\infty\gamma^k, \qquad k\ge 1,
\end{equation}
by the induction.
Let $B_i$ stand for $B(x, \rho^ir)$ and $\tau_i$ for $\tau_{ B(x, \rho^ir)}$.  Let
$$
m_i=\inf_{B_i}f, \qquad M_i=\sup_{B_i}f.
$$
Suppose $M_i-m_i\le \|f\|_\infty\gamma^i$ for all $i\le k$; we want to show that
\begin{equation}\label{e:4.4}
M_{k+1}-m_{k+1}\le \|f\|_\infty \gamma^{k+1}.
\end{equation}
Note that $m_k\le f\le M_k$ on $B_{k+1}$. Let
$$
A'=\{z\in B_{k+1}: f(z)\le \frac{m_k+M_k}2\}.
$$
We may assume $|A'|/|B_{k+1}|\ge 
1/2$, otherwise we look at 
$ \|f\|_\infty-f$ instead of $f$. 
Let $A$ be a compact subset of $A'$ with $|A|/|B_{k+1}|\ge 
1/3$. Let $\epsilon>0$ and choose  $y, z\in B_{k+1}$ with $f(y)
 \ge M_{k+1}-\epsilon$ and $f(z)\le m_{k+1}+\epsilon$.  Then
\begin{align*}
f(y)-f(z)=&\E_y[f(Y^\kappa_{T_A})-f(z); T_A<\tau_k]\\
&+\E_y[f(Y^\kappa_{\tau_k})-f(z); \tau_k<T_A, Y^\kappa_{\tau_k}\in B_{k-1}]\\
&+\sum^\infty_{i=1}\E_y[f(Y^\kappa_{\tau_k})-f(z); \tau_k<T_A, Y^\kappa_{\tau_k}\in B_{k-i-1}\setminus B_{k-i}]=:I+II+III.
\end{align*}
By the choice of $A$, 
$$
I \le \left(\frac{m_k+M_k}2-m_k\right)\P_y(T_A<\tau_k)=\frac12(M_k-m_k)\P_y(T_A<\tau_k)
$$
and, clearly, 
$$
II \le (M_k-m_k)\P_y(\tau_k<T_A)=(M_k-m_k)(1-\P_y(T_A<\tau_k)).
$$
By
the induction hypothesis,  
\eqref{e:4.2} and the fact that  $$\rho \le ({\gamma}/2)^{1/(2\delta_1)}
\wedge (c_1\gamma^2/(8c_3))^{1/(2\delta_1)},
$$
we have that 
\begin{align*}
& III \le \sum^\infty_{i=1}(M_{k-i-1}-m_{k-i-1})\P_y(Y^\kappa_{\tau_k}\in D \setminus B_{k-i})
\le \sum^\infty_{i=1}
c_3\|f\|_\infty\gamma^{k-i-1}\frac{(\rho^kr)^{2\delta_1}}{(\rho^{k-i}r)^{2\delta_1}}\\
&=
c_3\|f\|_\infty \gamma^{k-1}\sum^\infty_{i=1}(\rho^{2\delta_1}/\gamma)^i\le 2c_3\|f\|_\infty  \gamma^{k-2}\rho^{2\delta_1}\le \frac{c_1}4\|f\|_\infty \gamma^k.
\end{align*}
Therefore, by \eqref{e:4.1} and the fact that $\rho \le 1/3$, we have
\begin{align*}
f(y)-f(z)\le&\frac12
(M_k-m_k)\P_y(T_A<\tau_k)+(M_k-m_k)(1-\P_y(T_A<\tau_k))+\frac{c_1}4\|f\|_\infty \gamma^k\\
\le &(M_k-m_k)\left(1-\frac12\P_y(T_A<\tau_k)\right)+\frac{c_1}4\|f\|_\infty \gamma^k\\
\le &(M_k-m_k)\left(1-\frac12\P_y(T_A<\tau_{B(x, 3\rho^{k+1}r)})\right)+\frac{c_1}4\|f\|_\infty \gamma^k\\
\le &\|f\|_\infty\gamma^k\left(1-\frac{c_1}2\right)+\frac{c_1}4\|f\|_\infty\gamma^k
=\|f\|_\infty\gamma^k\left(1-\frac{c_1}4\right)=\|f\|_\infty\gamma^{k+1}.
\end{align*}
Hence
$$
 M_{k+1}-m_{k+1}\le f(y)-f(z) +2\epsilon\le \|f\|_\infty\gamma^{k+1}+2\epsilon.
$$
Since $\epsilon$ can be arbitrarily small,  \eqref{e:4.4} holds and hence 
\eqref{e:4.3} holds.

If $x, y\in B(x_0, \rho r/2)$, let $k$ be the smallest natural number with $|x-y|/r\le \rho^k$.  Then
$$
\log\frac{|x-y|}{r}\ge (k+1)\log\rho,
$$
$y\in B(x, \rho^kr)$, and
\begin{align*}
&|f(y)-f(x)|\le 
\|f\|_\infty \gamma^k=\|f\|_\infty e^{k\log\gamma}\\
&\le c_4\|f\|_\infty e^{\log\left(\frac{|x-y|}r\right)(\log\gamma/\log\rho)}
=c_4\|f\|_\infty \left(\frac{|x-y|}r\right)^{\log\gamma/\log\rho}.
\end{align*} 
If $x, y\in B(x_0, r/2) \setminus B(x_0, \rho r/2)$, then clearly $|f(x)-f(y)|\le 2\|f\|_\infty 
\le 
c_5 \|f\|_\infty \left(\frac{|x-y|}{r}\right)^\beta$.
\qed


\section{Existence of Green function }\label{s:exist-gf}
In the first part of this section, we assume that \textbf{(H1), }\textbf{(H2-w)}-\textbf{(H3-w)} hold and show that the process $Y^{\kappa}$ admits a Green function. 
Then we will assume additionally that \textbf{(H4-w)}-\textbf{(H5-w)} hold,  so that 
Proposition  \ref{t:uhp} holds.  Using Proposition \ref{t:uhp}, 
we will show that the Green function is finite off the diagonal.

First we assume that \textbf{(H1)},  \textbf{(H2-w)}-\textbf{(H3-w)}  hold.
Recall that  $\zeta$ is the  lifetime of $Y^{\kappa}$.
Let $f:D\to [0,\infty)$ be a Borel function and $\lambda \ge 0$. 
The $\lambda$-potential of $f$ is defined by
$$
G_{\lambda} f(x):=\E_x\int_0^\zeta e^{-\lambda t}f(
Y^{\kappa}_t
)\, dt\, , \quad x\in D.
$$
When $\lambda=0$, we write $Gf$ instead of $G_0f$ and call $Gf$ the Green potential of $f$. If $g:D\to [0,\infty)$  is another Borel function, then by the symmetry of $
Y^\kappa$ we have that 
\begin{equation}\label{e:symmetry-of-G}
\int_{D}G_{\lambda} f(x)g(x)\, dx = \int_{D}f(x) G_{\lambda} g(x)\, dx \, .
\end{equation}
For $A\in \BB(D)$, we let $G_{\lambda}(x,A):=G_{\lambda} \ind_A (x)$ be the $\lambda$-occupation measure of $A$.

Let $U$ be a relatively compact 
open subset of $D$. 
For $\gamma>0$, let $J_\gamma$ be the jump kernel defined 
in the proof of Lemma \ref{l:YDkappaU} and let $Z$ be the pure jump conservative process with jump kernel $J_\gamma$.  In the proof of Lemma \ref{l:YDkappaU}  we have shown that, when $\gamma$ is small enough, the function $\wt\kappa$ defined in \eqref{e:wt-kappa} is non-negative and the semigroup $(Q_t^U)_{t\ge 0}$ of $Y^{\kappa, U}$ is given by 
$$
Q_t^U f(x)=\E_x[\exp(-A_t)f(Z_t^U)]\, , \quad t>0, \ x\in U,
$$
where $A_t:=\int_0^t\wt{\kappa}(Z_s^U)\, ds$.
Moreover,  $Q_t^U$ has a  transition density $q^U(t,x,y)$ (with respect to the Lebesgue measure) which is symmetric in $x$ and $y$, and such that for all $y\in U$, $(t,x)\mapsto q^U(t,x,y)$ is continuous. 

Let $G_{\lambda}^U f(x):=\int_0^{\infty}e^{-\lambda t}Q_t^U f(x)\, dt=\E_x \int_0^{\tau_U}e^{-\lambda t} f(
Y^{\kappa}_t
)\, dt$ denote the $\lambda$-potential of $Y^U$ and  $G_{\lambda}^U(x,y):=\int_0^{\infty}e^{-\lambda t}q^U(t,x,y)\, dt$ the $\lambda$-potential density of $Y^U$. 
We will write $G^U$ for $G_{0}^U$ for simplicity.
Then 
$G_{\lambda}^U(x,\cdot)$ is the density of the $\lambda$-occupation measure. In particular this shows that $G_{\lambda}^U(x,\cdot)$ is absolutely continuous with respect to the Lebesgue measure.
Moreover, since $x\mapsto q^U(t,x,y)$ is continuous, we see that $x\mapsto G^U_{\lambda}(x,y)$ is lower semi-continuous. By Fatou's lemma this implies that $G_{\lambda}^U f$ is also lower semi-continuous.

Let $(U_n)_{n\ge 1}$ be a sequence of bounded 
open sets such that $U_n\subset \overline{U_n}\subset U_{n+1}$ and $\cup_{n\ge 1}U_n=D$. 
For any Borel $f: D\to [0,\infty)$, it holds that

\begin{equation}\label{e:Gf-increasing-limit}
G_{\lambda}f(x)=\E_x \int_0^{\zeta}e^{-\lambda t}f(
Y^{\kappa}_t
)\, dt  =\, \uparrow\!\!\! \lim_{n\to \infty}\E_x \int_0^{\tau_{U_n}}e^{-\lambda t}f(
Y^{\kappa}_t)\, dt =\, \uparrow\!\!\!\lim_{n\to \infty}G_{\lambda}^{U_n}f(x)\, ,
\end{equation}
 where $ \uparrow\!\!\!\ \lim $ denotes an increasing limit. 

In particular, if 
$A\in \BB(D)$ is of Lebesgue measure zero, then for every $x\in D$,
$$
G_{\lambda}(x,A)=\lim_{n\to \infty}G_{\lambda}^{U_n}(x,A)=\lim_{n\to \infty}G_{\lambda}^{U_n}(x,A\cap U_n)=0\, .
$$
Thus,  $G_{\lambda}(x,\cdot)$ is absolutely continuous with respect to the Lebesgue measure 
for each $\lambda\ge 0$ and $x\in D$. Together with \eqref{e:symmetry-of-G} 
this shows that the conditions of \cite[VI Theorem (1.4)]{BG} are 
satisfied, which  implies 
that the resolvent $(G_{\lambda})_{\lambda>0}$ is self dual. In particular, 
see \cite[pp.256--257]{BG}, 
there exists a symmetric function $G(x,y)$ excessive in both variables such that
$$
Gf(x)=
\int_D
G(x,y)f(y)\, dy\, ,\quad x\in 
D.
$$
We recall, see \cite[II, Definition (2.1)]{BG}, that a measurable
function $f:
D\to [0,\infty]$ 
is $\lambda$-excessive, $\lambda \ge 0$, with respect to the process $Y^{\kappa}$ if for every $t\ge 0$ it holds that 
$\E_x[e^{-\lambda t}Y^{\kappa}_t]\le f(x)$ and $\lim_{t\to 0}\E_x[e^{-\lambda t}Y^{\kappa}_t] =f(x)$, 
for every $x\in D$.
$0$-excessive functions are simply called excessive functions.

We note that the process $Y^{\kappa}$ need not be transient. If it is transient, then it follows that $G(x,y)<\infty$ for a.e.~$y\in D$. In the 
following 
lemma we show transience under the additional assumption that $\kappa$ is strictly positive.

\begin{lemma}\label{l:transience}
Suppose that $\kappa(x)>0$ for every $x\in D$. Then the process $Y^{\kappa}$ is transient in the sense that there exists $f:D\to (0, \infty)$ such that $Gf<\infty$. More precisely, $G\kappa\le 1$.
\end{lemma}
\pf Let $(Q_t)_{t\ge 0}$ denote the semigroup of $Y^{\kappa}$. For any $A\in \BB(D)$, we use \cite[(4.5.6)]{FOT} with $h=\ind_A$, $f=1$, and let $t\to \infty$ to obtain
$$
 \int_A \P_x(\zeta <\infty)\, dx \ge \int_A \P_x(Y^{\kappa}_{\zeta-} \in D,  \zeta <\infty)\, dx 
=\int_0^{\infty}
\int_D
\kappa(x) Q_s \ind_A(x)\, dx \, dt.
$$
This can be rewritten as
$$
\int_A
 \P_x(\zeta<\infty)
\, dx \ge \int_{D} \kappa(x) G\ind_A (x)\, dx =\int_A G\kappa(x)\, dx.
$$
Since this inequality holds for every $A\in \BB(D)$, we conclude that 
$
\P_x( \zeta<\infty)
 \ge G\kappa(x)$ for a.e.~$x\in D$. 
Both functions $x\mapsto \P_x(\zeta<\infty)$ and $G\kappa$ are excessive. Since $G(x, \cdot)$ is absolutely continuous with respect to the Lebesgue measure 
(i.e.,~Hypotesis (L) holds, see \cite[p.112]{CW}), 
by  \cite[Proposition 9, p.113]{CW}, 
we conclude that  $G\kappa(x) \le \P_x(\zeta <\infty)\le 1$ for all $x\in D$. 
\qed

From now on we assume that $Y^{\kappa}$ is transient so that $G(x,y)$ is not identically infinite.
Note that it follows from \eqref{e:Gf-increasing-limit} that, for every non-negative Borel $f$, $G_{\lambda}f$ is lower semi-continuous, as an increasing limit of lower semi-continuous functions. Since every $\lambda$-excessive function is an increasing limit of $\lambda$-potentials see  \cite[II Proposition (2.6)]{BG}), 
we conclude that all $\lambda$-excessive functions of $
Y^{\kappa}$ are lower semi-continuous. In particular, for every $y\in 
D$, $G_{\lambda}(\cdot, y)$ is lower semi-continuous. Since $G(\cdot, y)$ is the increasing limit of $G_{\lambda}(\cdot, y)$ as $\lambda\to 0$, we see that $G(\cdot, y)$ is also lower semi-continuous.

Fix an open set $B$ in $D$ and $x \in D$.
Let $f$ be a non-negative Borel function on $D$. By Hunt's switching identity, \cite[VI, Theorem (1.16)]{BG}, 
$$
\E_x [G f(
Y^{\kappa}_{\tau_B})]=
\int_D
\E_x [G(
Y^{\kappa}_{\tau_B},y)] f(y)\, dy =
\int_D
\E_y [G(x,Y^{\kappa}
_{\tau_B})] f(y)\, dy. 
$$
Suppose, further, that $f=0$ on $B$. Then by the strong Markov property, 
\cite[I, Definition (8.1)]{BG},
$$
\int_D
G(x,y) f(y)\, dy
= \E_x \int_{\tau_B}^{\infty} f(
Y^{\kappa}_t)\, dt 
= \E_x[G f(
Y^{\kappa}_{\tau_B})]=
\int_{D\setminus B} 
\E_y [G(x,
Y^{\kappa}_{\tau_B})] f(y)\, dy\, ,$$
and hence $G(x,y)=\E_y [G(x,
Y^{\kappa}_{\tau_B})]$ for a.e. $y\in D \setminus B$. Since both sides are excessive (and thus excessive for the killed process $
Y^{\kappa, D\setminus B}
$), equality holds for every $y\in D \setminus B$. By using Hunt's switching identity one more time, we arrive at
$$
G(x,y)=\E_x [G (
Y^{\kappa}_{\tau_B},y)]\, ,\quad \text{for all } x\in D,\ y\in D\setminus B\, .
$$
In particular, if $y\in D\setminus B$ is fixed, then the above equality says that $x\mapsto G(x,y)$ is regular harmonic in $B $ with respect to $
Y^{\kappa}$. By symmetry, $y\mapsto G(x,y)$ is regular harmonic in $B $ as well.

Now we assume additionally that  \textbf{(H4-w)}-\textbf{(H5-w)} hold.
 By using Proposition \ref{t:uhp} \bk
we conclude that $G(x,y)<\infty$ for all $y\in D\setminus \{x\}$. 
This proves  the following 
 result about 
the existence of the Green function.

\begin{prop}\label{p:existenceGF}
Suppose that \textbf{(H1)},  \textbf{(H2-w)}-\textbf{(H5-w)}, 
 \eqref{e:assumption-on-j}-\eqref{e:H-infty} and \eqref{e:bound-on-kappa0} hold.
Assume that $Y^{\kappa}$ is transient.
Then there exists a symmetric function $G:D\times D\to [0,\infty]$ which
is lower semi-continuous in each variable and finite outside the diagonal such that for every non-negative Borel $f$,
$$
Gf(x)=
\int_D
G(x,y)f(y)\, dy\, .
$$
Moreover, $G(x, \cdot)$ is harmonic with respect to $Y^{\kappa}$ in $D\setminus \{x\}$
and regular harmonic with respect to $
Y^{\kappa}$ in $D\setminus B(x, \epsilon)$ for any $\epsilon>0$.
\end{prop}

We now prove the continuity of Green function under an additional assumption. 
The proof of the next proposition is similar to the corresponding part of the proof of \cite[Theorem 1.1]{KSV21G}.
\begin{prop}\label{p:conti_GF}
Suppose that 
\textbf{(H1)}-\textbf{(H4)},  \textbf{(H5-w)},
\eqref{e:assumption-on-j}-\eqref{e:H-infty} and \eqref{e:bound-on-kappa} hold. 
Assume that $Y^{\kappa}$ is transient and that the Green function $G:D\times D\to [0,\infty]$ of $Y^{\kappa}$ 
satisfies that for any $x \in D$ and $r>0$ \begin{align}\label{e:Gassump}
\sup_{z\in D \setminus B(x, r)}G(x,z) <\infty.
\end{align}
Then $G(x, \cdot)$ is continuous in $D\setminus \{x\}$.
\end{prop}
\pf
We fix $x_0, y_0 \in D$ , $x_0\neq y_0$, 
and choose a positive $a$ small enough  so that $B(x_0, 4a) \cap B(y_0, 4a)= \emptyset$ and
$B(x_0, 4a) \cup B(y_0, 4a) \subset D$.

We first note that 
for $(z, w) \in B(x_0, 2a) \times B(y_0, 2a)$, 
 $\delta_D(z)\wedge \delta_D(w)\ge 2a = \frac{2a}{|x_0-y_0|+4a}(|x_0-y_0|+4a) \ge \frac{2a}{|x_0-y_0|+4a}|w-z|$. 
 Thus, by \textbf{(H2)},
\begin{align}
\sup_{(z, w) \in B(x_0, 2a) \times B(y_0, 2a)} J(z, w) \le 
c_0 
  \sup_{(z, w) \in B(x_0, 2a) \times B(y_0, 2a)} j(z, w)
\le \frac{c_1} {a^{d} \Psi(a)}.
\end{align}
We recall that by Proposition \ref{p:exit-time-estimate}(b), $\E_y\tau_{B(x_0, 2a)}\le c_2 \Psi(a)$ for all $y\in B(x_0,a)$. 
Let $N\ge 1/a$. 
In the paragraph after the proof of Lemma \ref{l:transience}, we have seen that for any non-negative Borel function $f$ and 
$\lambda\ge 0$, $G_\lambda f$ is lower semi-continuous. Thus by \cite[Theorem 2, p.126]{CW}, $G$ is locally integrable  in each variable.
Using \eqref{e:levys} in the second line and the local integrability of $G$ in the fourth, we have for every $y\in B(x_0,a)$, 
\begin{align*}
&\E_y\left[G(Y^\kappa_{\tau_{B(x_0, 2a)}}, y_0); Y^\kappa_{\tau_{B(x_0, 2a)}} \in B(y_0, 1/N) \right]\\
& = \E_y \left(\int_0^{\tau_{B(x_0, 2a)}} 
 \int_{B(y_0, 1/N)}G(w, y_0)J(Y^\kappa_s,w)dw\, ds\right)\\
& \le 
\left(\sup_{y \in B(x_0, a)}
\E_y\tau_{B(x_0, 2a)}\right) \left(\sup_{
z \in B(x_0, 2a)}
\int_{B(y_0, 1/N)} J(z, w)G(w, y_0)dw \right)\\
 & \le c_1c_2 a^{-d} \int_{B(y_0, 1/N)} G(w, y_0)dw <\infty.
\end{align*}
Given $\epsilon>0$,  choose $N$ large enough such that  $c_1c_2 a^{-d} \int_{B(y_0, 1/N)} G(w, y_0)dw<\epsilon/4$, so
$$
\sup_{y \in B(x_0, a)}\E_y
\left[G(Y^\kappa_{\tau_{B(x_0, 2a)}}, y_0); Y^\kappa_{\tau_{B(x_0, 2a)}} \in B(y_0, 1/N) \right]<\epsilon/4.
$$
The function 
$y\mapsto h(y) 
:= \E_y\left[G(Y^\kappa_{\tau_{B(x_0, 2a)}}, y_0); Y^\kappa_{\tau_{B(x_0, 2a)}}\in D \setminus B(y_0, 1/N)\right]$ is 
 harmonic on $ B(x_0, a)$, and  by \eqref{e:Gassump} 
it is bounded function on $D$. Thus,  by  Theorem \ref{t:holder},
 it is continuous.  Choose a $\delta 
\in (0, a)$ such that $|h(y)-h(x_0)|<\epsilon/2$ for all 
$y\in B(x_0, \delta)$ ,
We now see that for all  $y\in B(x_0, \delta)$,
\begin{align*}
&|G(y, y_0)-G(x_0, y_0)| \\
& \le |h(y)-h(x_0)|+2
\sup_{y \in B(x_0, a)}\E_y
\left[G(Y^\kappa_{\tau_{B(x_0, 2a)}}, y_0); Y^\kappa_{\tau_{B(x_0, 2a)}} \in B(y_0, 1/N) \right] <\epsilon.
\end{align*}
\qed


\section{Examples}\label{s:examples} 
In this section we give two families of examples of jump kernels that satisfy hypotheses 
\textbf{(H1)}-\textbf{(H5)}.

\subsection{Trace processes and resurrected kernels} \label{ss:E1}
Let $X=(X_t, \P_x)$ be a L\'evy process with  L\'evy measure $j(|x|)dx$. For the moment we do not assume that \eqref{e:assumption-on-j} and \eqref{e:H-infty} hold. 
Let $D$ be a proper open set in $\R^d$ 
 such that $U:=\overline D^c$ is non-empty. 
We denote the jump kernel of $X$ as $j(x,y)=j(|x-y|)$. 
Let
$$
A_t:=\int_0^t \1_{(X_s\in D)}\, ds
$$
and let $\tau_t:=\inf\{s>0: \, A_s>t\}$ 
be its right-continuous inverse. The process 
$Y=(Y_t)_{t\ge 0}$ 
defined by $Y_t:=X_{\tau_t}$ is a Hunt process with state space $D$.  The process  $Y$ is called the trace process of $X$ on $D$ (it is also called
the path-censored process in some literature, for instance, \cite{KPW14}).
Here is another way to describe 
the part of the process $Y$ until its first hitting time to the boundary $\partial D$:
Let $x=X_{\tau_D-}\in D$ be the position from which $X$ jumps out of $D$, and let 
$z=X_{\tau_D}\in U$ be the position where $X$ lands at the exit from $D$. The distribution of the returning position of $X$ to $D$ is given by the Poisson kernel of $X$ 
with respect to $U$:
$$
P_U(z, A)=\int_A \int_U G_U(z,w)j(w,y)\,  dw\, dy, \qquad A\in \BB(D).
$$
Here $G_U(z,w)$, $z,w\in U$, denotes the Green function of the process $X$ killed upon exiting $U$.
This implies that when $X$ jumps out of $D$ from the point $x$, we continue the process 
by resurrecting it in $ A\in \BB(D)$ according to the kernel
$$
q(x,A)=\int_{U}j(x,z)
P_U(z, A)
\, dz, \quad x\in D,
$$
which has  density
$$
q(x,y)=\int_U \int_U j(x,z)G_U(w,z)j(z,y)\, dz\, 
dw, \quad x,y\in D.
$$
We call $q(x,y)$ the \emph{resurrection kernel}. Since the Green function $G_U$ is symmetric, it immediately follows that $q(x,y)=q(y,x)$ for all $x, y\in D$. This shows that 
the part of the process $Y$ until its first hitting of the boundary 
can be regarded as a resurrected process. The jump kernel of this process is symmetric and is given by $J(x,y)=j(x,y)+q(x,y)$, $x,y\in D$.

This example can be modified in the following way. For each $z\in U$, let $p(z, y)$ be a subprobability density on $D$. Instead of returning the process $X$ to $D$ by using the Poisson kernel $P_U(z,A)$, we 
may
use the kernel $p(z,A)=\int_A p(z,y)dy$, $A\in \BB(D)$. We call this kernel the \emph{return} kernel. Define
\begin{equation}\label{e:defq}
q(x,y):=\int_U j(x,z)p(z,y)\, 
dz.
\end{equation}
We assume that $p(z,y)$ satisfies the following two properties: (1) It is such that $q$ is symmetric, that is, $q(x,y)=q(y,x)$; (2) 
There exists $c_1\ge 1$ such that for all $y_0\in D$ and $r>0$ with $B(y_0, 2r)\subset D$,
\begin{equation}\label{e:dfp}
c_1^{-1}p(w,y_1)\le p(w, y_2) \le c_1 p(w,y_1) \quad \text{for all } w\in D^c \text{ and all } y_1,y_2\in B(y_0,r).
\end{equation}
Clearly, both properties are true for the trace process (that is, for the Poisson kernel $P_U(w,y)$). 
We note that the first property is quite delicate. Still, many examples of return kernels for which $q$ is symmetric are given in \cite{KSV22b} for the case of a half-space. 
Both properties 
can be checked in concrete examples of return kernels. One such example is given by
$$
p(w,y):=\frac{j(w,y)}{\int_D j(w,z)\, dz}
$$
studied in 
\cite{DRV, Von21} 
in the context of the Neumann boundary problem. Note that it follows from \eqref{e:j12} below that this $p$ satisfies 
both properties above.
 
 For $x,y\in D$, $x\neq y$, define
\begin{align}
\label{e:JRut}
J(x,y):=j(x,y)+q(x,y)=j(x,y)\left(1+\frac{q(x,y)}{j(x,y)}\right)=:j(x,y)\sB(x,y), 
\end{align}
where
\begin{align}\label{e:RsB}
\sB(x,y)=
\begin{cases}
1+\frac{q(x,y)}{j(x,y)} \quad &\text{ when } x\not=y;\\
1 \quad &\text{ when } x=y.
\end{cases}
\end{align}
In the remaining part of this subsection we show that $J(x,y)$ (that is, $\sB(x,y)$) satisfies hypotheses 
\textbf{(H1)}-\textbf{(H5)}.

\medskip
Since we have assumed that the return kernel $p(z,y)$ is such that $q$ is symmetric, it is immediate that $\sB(x,y)=\sB(y, x)$ for all $x,y\in D$.
Hence \textbf{(H1)}  holds.

\medskip
Fix  $\epsilon \in (0,1)$, $x_0\in \R^d$ and $r>0$.
Then, for all $w \in \R^d\setminus B(x_0, (1+\epsilon)r)$ and $x_1,x_2\in B(x_0,r)$,
$$|x_1-w| \le  |x_2-w| +|x_1-x_2|\le |x_2-w| +2r
\le(1+(2/\epsilon))|x_2-w|.
$$
Thus $j(x_1,w) \asymp j(x_2,w)$ for $w \in \R^d \setminus B(x_0, (1+\epsilon)r)$ and $x_1,x_2\in B(x_0,r)$. 
In particular, if $x_0\in D$ and $B(x_0, (1+\epsilon)r)\subset D$, then 
$$
j(x_1,z) \asymp j(x_2,z)\quad \text{for all }x_1,x_2\in B(x_0,r), \, \,z\in D\setminus B(x_0, (1+\epsilon)r)
$$
and
\begin{align}\label{e:j12}
j(x_1,w) \asymp j(x_2,w)\quad \text{for all }x_1,x_2\in B(x_0,r), \, \,w\in D^c.
\end{align}
Therefore,  $\text{for all }x_1,x_2\in B(x_0,r), \, \,z\in D\setminus B(x_0, (1+\epsilon)r),$
\begin{align*}
\sB(x_1,z)&=
1+\frac{1}{j(x_1,z)}\int_{U}j(x_1,w)p(w,z)\, 
dw \\
&\asymp
1+\frac{1}{j(x_2,z)}\int_{U}j(x_2,w)p(w,z)\, 
dw=\sB(x_2,z).
\end{align*}
Hence  \textbf{(H5)} holds. 

\medskip
To check \textbf{(H2)} and \textbf{(H4)}, we will use the following two lemmas for $q$.
\begin{lemma}\label{l:qcom}
For every $\epsilon \in (0,1)$, there exists 
$C_{25}=C_{25}(\epsilon)\ge 1$ such that 
for all $x_0, y_0\in D$ and $r>0$, 
$$
C_{25}^{-1} q(x_0, y_0) \le q(x,y) \le C_{25}
q(x_0, y_0), \quad  (x,y) \in B(x_0, (1-\epsilon)\delta_D(x_0)) \times B(y_0, (1-\epsilon)\delta_D(y_0)) . 
$$
\end{lemma}
\pf 
The lemma  follows from \eqref{e:dfp} and \eqref{e:j12}. \qed

\begin{lemma}\label{l:qubd}
There exists $C_{26}>1$ such that 
$$
q(x, y) \le 
C_{26} 
\left(\frac{1}{\delta_D(y)^d\Psi(\delta_D(x))} \wedge \frac{1}{\delta_D(x)^d\Psi(\delta_D(y))} \right)
\quad \text{ for all }x, y\in D.$$
\end{lemma}
\pf 
By Lemma \ref{l:qcom},
$$q(x, y) \asymp q(x, u)=\int_U 
j(x,w)p(w,u)dw\quad \text{for all } u \in B(y, \delta_D(y)/2).
$$ 
Thus, by \eqref{e:ww1},
\begin{align*}
q(x, y) 
&\le \frac{c_1}{\delta_D(y)^d} \int_{B(y, \delta_D(y)/2)}q(x, u)du=\frac{c_1}{\delta_D(y)^d}\int_U j(x,w) \left(\int_{B(y, \delta_D(y)/2)}p(w,u)du \right)dw\\
&\le \frac{c_1}{\delta_D(y)^d} \int_{B(x, \delta_D(x))^c} j(x,w)dw \le \frac{c_2}{\delta_D(y)^d\Psi(\delta_D(x))}. 
\end{align*}
The lemma now follows from this and the symmetry of $q$.
\qed

 Applying Lemma \ref{l:qubd} to \eqref{e:RsB} and using \eqref{e:assumption-on-j}, 
 we get
\begin{align*}
\sB(x,y)-1&\le
\frac{c_1}{(\delta_D(x)\wedge \delta_D(y))^d\Psi(\delta_D(x)\wedge \delta_D(y))j(|x-y|)}\\
&\asymp  
\frac{\Psi(|x-y|)}{\Psi(\delta_D(x)\wedge \delta_D(y))} \left(\frac{|x-y|}{\delta_D(x)\wedge \delta_D(y)}  \right)^d
 \le c_2 \left(\frac{|x-y|}{\delta_D(x)\wedge \delta_D(y)}  \right)^{d+2\delta_1}.
\end{align*}
This proves  that  both \textbf{(H2)} and \textbf{(H4)} hold.

\medskip

To check \textbf{(H3)}, it suffices to show \eqref{e:H3} for $a\in (0,1/2]$. 
Let $j(x,dz):=j(x,z)dz$. Then  for any $x \in D$, $ j(x,dz)$ is a finite measure on $D^c$ 
such that, by \eqref{e:ww1},
\begin{align}\label{e:Jw}
j(x,D^c) \le  \frac{c_3}{\Psi(\delta_D(x))}  
\qquad x \in D,
\end{align}
for some $c_3>0$.
By \eqref{e:ww1}, \eqref{e:defq}, \eqref{e:JRut}, \eqref{e:Jw}, and the fact that $p$ is a subprobability kernel, 
for $a\in (0,1/2]$, 
\begin{align*}
&\int_{D, |x-y|>a\delta_D(x)} J(x,y)dy
 \le \int_{|x-y|>a\delta_D(x)} j(x,y)dy+\int_{D} \int_{D^c} j(x,dw) p(w,y)\, dy \nn \\
& \le 
c_4 \Psi(a\delta_D(x))^{-1} + \int_{D^c}\left(\int_{D} p(w,y) dy\right) j(x,dw)  \le c_5(a)\Psi(\delta_D(x))^{-1}.
\end{align*}
This proves that \textbf{(H3)} holds.

\begin{remark}{\rm
Suppose that for every $x\in D$, $\tilde{j}(x,dz)$ is a kernel on $D^c$ satisfying \eqref{e:Jw}. For $x,y\in D$, let $\tilde{q}(x,y):=\tilde{j}(x,dz)p(z,y)$. Then $J(x,y):=j(x,y)+\tilde{q}(x,y)$ also satisfies \textbf{(H3)}.
}
\end{remark}

\subsection{Examples of $\sB(x,y)$  satisfying \textbf{(H3)} which may blow up  at the boundary.}\label{ss:E2}

Let $D\subset \R^d$ be a proper open subset of $\R^d$, $J(x,y)=j(|x-y|)\sB(x,y)$, where $j$ satisfies \eqref{e:assumption-on-j} with 
$\Psi$ satisfying \eqref{e:H-infty}.

We first record an estimate of  $j(|x-y|)$
 in case $|x-y| > a \delta_D(x)$. By \eqref{e:H-infty},
$$
\Psi(|x-y|)^{-1}\le a_1^{-1}\left(\frac{a \delta_D(x)}{|x-y|}\right)^{2\delta_1} \Psi(a \delta_D(x))^{-1}
$$
and 
$$
\Psi(a \delta_D(x))^{-1}\le a_2 a^{-2\delta_2}\Psi(\delta_D(x))^{-1}.
$$
This, together with \eqref{e:assumption-on-j}, implies that there exists $c_1=c_1(a,\delta_1, \delta_2)>0$ such that
\begin{align}
 \label{e:jj}
j(|x-y|)\le c_1 \Psi(\delta_D(x))^{-1} \delta_D(x)^{2\delta_1} |x-y|^{-d-2\delta_1}, \quad |y-x|> a \delta_D(x).
\end{align}

\begin{lemma}\label{l:estimate-for-halfspace-H3}
Suppose that 
$D$ is a proper open subset of $\R^d$
and let $L_s:=\{y\in D: \delta_D(y)=s\}$. 
When $d\ge 2$ we  assume that  
 there exists $A_1>0$ such that
\begin{equation}\label{e:Hcond}
\HH_{d-1}(L_s \cap B(z,  R))\le A_1 R^{d-1}, \quad z\in D,  s>0, R>0,
\end{equation}
where  $
\HH_{d-1}$ is the $(d-1)$-dimensional 
Hausdorff  measure.
Assume
that there exists $\beta_2\in [0, 1\wedge (2\delta_1))$ 
with the property that for all $a\in (0,1/2]$, 
 there exists $A_2(a)>0$ such that 
\begin{equation}\label{e:assumption-B-example-2}
\sB(x,y)\le 
A_2(a) \frac{|x-y|^{2\beta_2}}{\delta_D(x)^{\beta_2}\delta_D(y)^{\beta_2}}, 
\quad \text{ for }|x-y|> a\delta_D(x).
\end{equation} 
 Then for every  $a\in (0,1/2]$, there exists $C_{27}(a)>0$ such that
$$
\int_{|x-y|> a \delta_D(x)} J(x,y)\, dy \le 
C_{27}(a) \Psi(\delta_D(x))^{-1}.
$$
\end{lemma}
\pf
We give the proof for $d\ge 2$.
It follows from \cite[Theorem 6.3.3 (vi) and (vii), p. 285]{DZ} that the function 
$x \mapsto \delta_D(x)$ is Lipschitz on $D$ and $|\nabla \delta_D(x)|=1$
a.e. $ x\in D$.
  Thus,  the following coarea formula
is valid (see \cite[Theorem 3.2.3 (2)]{F}: For any $g\in L^1(D)$,
\begin{equation}\label{e:coarea}
\int_Dg(y)dy=\int^\infty_0\int_{
L_s}g(y)
\HH_{d-1}(dy)ds.
\end{equation}

It follows from \eqref{e:jj} and \eqref{e:assumption-B-example-2} that
\begin{equation}\label{e:basic-estimate-for-H3}
\int_{|x-y|>a \delta_D(x)}J(x,y)\, dy \le c_1 
\frac{ \delta_D(x)^{2\delta_1-\beta_2}}{\Psi(\delta_D(x))}
\int_{|x-y|> a \delta_D(x)}
\frac{\delta_D(y)^{-\beta_2} }{|x-y|^{d+2\delta_1 -2\beta_2}}\, dy.
\end{equation}

We split the integral into two parts:
\begin{align*}
& \int_{|x-y|> a \delta_D(x)}
\frac{\delta_D(y)^{-\beta_2} }{|x-y|^{d+2\delta_1 -2\beta_2}}\, dy\\
& =\int_{|x-y|>a\delta_D(x), \delta_D(y) \le (1+a)\delta_D(x)}+\int_{\delta_D(y)>(1+a)\delta_D(x)}\\
& =:I_1+I_2. 
\end{align*}
Here we used that if 
$\delta_D(y) \ge(1+a)\delta_D(x)$,
 then $|x-y|> a \delta_D(x)$.

Using  \eqref{e:Hcond}, \eqref{e:coarea}, 
the assumption $\beta_2<1\wedge (2\delta_1)$
in the last line, we have
\begin{align*}
&I_1=
\int_D{\bf 1}_{B(x, a\delta_D(x))^c}(y){\bf 1}_{ \delta_D(y)\le (1+a)\delta_D(x)\}}
\frac{\delta_D(y)^{-\beta_2}dy}{|x-y|^{d+2\delta_1 -2\beta_2}}\\
&=
\sum^\infty_{n=0}
\int_D{\bf 1}_{\delta_D(y)\le (1+a)\delta_D(x)}{\bf 1}_{B(x, 2^{n+1}a\delta_D(x))\setminus B(x, 2^{n}a\delta_D(x)}(y)
\frac{\delta_D(y)^{-\beta_2}dy}{|x-y|^{d+2\delta_1 -2\beta_2}}
\\
&\le
\delta_D(x)^{-d-2\delta_1+2\beta_2}\sum^\infty_{n=0}2^{n(-d-2\delta_1 +2\beta_2)}
\int_D{\bf 1}_{\delta_D(y)\le (1+a)\delta_D(x)}{\bf 1}_{B(x, 2^{n+1}a\delta_D(x))}(y)\delta_D(y)^{-\beta_2}
dy\\
&
=
\delta_D(x)^{-d-2\delta_1+2\beta_2}\sum^\infty_{n=0}2^{n(-d-2\delta_1 +2\beta_2)}\int^{(1+a)\delta_D(x)}_0
\HH_{d-1}(L_s\cap B(x,2^{n+1}a\delta_D(x)))s^{-\beta_2}ds\\
&\le c_2\delta_D(x)^{-1-2\delta_1+2\beta_2}\sum^\infty_{n=0}
2^{n(-d-2\delta_1 +2\beta_2)}2^{(n+1)(d-1)}\int^{(1+a)\delta_D(x)}_0 s^{-\beta_2}ds\\
&\le c_3\delta_D(x)^{-2\delta_1 +\beta_2}\sum^\infty_{n=0}2^{n(
-2\delta_1
 +2\beta_2-1)}=
 c_4\delta_D(x)^{-2\delta_1 +\beta_2}.
\end{align*}

When $s>(1+a)\delta_D(x)$, we have $|x-y| \ge s-\delta_D(x) \ge (a/(1+a)) s$.
Using  \eqref{e:Hcond}, \eqref{e:coarea}, 
the facts that $2\beta_2<1+2\delta_1$ and $\beta_2<2\delta_1$ in the last line, we have
\begin{align*}
 I_2 =&
\int_{D}{\bf 1}_{\{\delta_D(y)>(1+a)\delta_D(x) \}}
\frac{\delta_D(y)^{-\beta_2} }{|x-y|^{d+2\delta_1 -2\beta_2}}\, dy\\
=&
\int_{(1+a)\delta_D(x)}^\infty
\int_{L_s}
|x-y|^{-d-2\delta_1 +2\beta_2} \, 
\HH_{d-1}(dy) s^{-\beta_2}ds\\
=&
\int_{(1+a)\delta_D(x)}^\infty
\int_{L_s, |x-y| \le 2s }
|x-y|^{-d-2\delta_1 +2\beta_2} 
\HH_{d-1}(dy) s^{-\beta_2}ds\\
&+\int_{(1+a)\delta_D(x)}^\infty
 \int_{L_s, |x-y| > 2s }
|x-y|^{-d-2\delta_1 +2\beta_2} 
\HH_{d-1}(dy) s^{-\beta_2}ds\\
\le &
\int_{(1+a)\delta_D(x)}^\infty
[(a/(1+a)) s]^{-d-2\delta_1 +2\beta_2}
\HH_{d-1}(L_s\cap B(x, 2s))  s^{-\beta_2}ds\\
&+\int_{(1+a)\delta_D(x)}^\infty \sum_{n=1}^\infty
 \int_{
 L_s, 2^{n+1}s\ge |x-y| > 2^{n}s } 
 |x-y|^{-d-2\delta_1 +2\beta_2}\, 
\HH_{d-1}(dy) s^{-\beta_2}ds\\
\le &c_5 [(a/(1+a))]^{-d-2\delta_1 +2\beta_2}
\int_{(1+a)\delta_D(x)}^\infty s^{-2\delta_1 +\beta_2-1}ds\\
&+\int_{(1+a)\delta_D(x)}^\infty \sum_{n=1}^\infty
(2^{n}s)^{-d-2\delta_1 +2\beta_2}\, 
\HH_{d-1}(L_s\cap B(x, 2^{n+1}s)) 
s^{-\beta_2}ds
\\
\le &
c_6 \delta_D(x)^{-2\delta_1 +\beta_2}
+c_7\sum_{n=1}^\infty
2^{(-2\delta_1 +2\beta_2-1)n}  \int_{(1+a)\delta_D(x)}^\infty s^{-1-2\delta_1 +\beta_2}ds
 \le c _8\delta_D(x)^{-2\delta_1 +\beta_2}.
\end{align*}
Combining the display above with \eqref{e:basic-estimate-for-H3}, we get the conclusion of the lemma.
\qed

In the remaining part of this subsection we impose the following conditions on $\sB(x,y)$ that 
is  used in \cite[Section 4]{KSV22b} in the case of the half-space. Suppose $0\le \beta_1\le \beta_2 < 1\wedge (2\delta_1)$. Let $\Phi$ be a positive function on 
 $(0, \infty)$ satisfying 
 $\Phi(t)  \equiv \Phi(2)>0$ on $[0, 2)$ and  the following weak scaling condition: 
There exist constants $b_1, b_2>0$ such that 
\begin{equation}\label{e:Phi-infty}
   b_1(R/r)^{\beta_1}\leq\frac{\Phi(R)}{\Phi(r)}\leq b_2(R/r)^{\beta_2},\quad  2\le r<R <\infty.
\end{equation}

Recall that $D$ is a proper open subset of $\R^d$.
Assume that $\sB(x,y)$ satisfies \textbf{(H1)}, \textbf{(H4)} and the following assumption: There exists 
$C_{28} \ge 1$ such that
\begin{equation}\label{e:B7}
C_{28}^{-1} \Phi\left(\frac{|x-y|^2}{\delta_D(x)\delta_D(y)}\right)\le \sB(x,y)
 \le C_{28}\Phi\left(\frac{|x-y|^2}{\delta_D(x)\delta_D(y)}\right) \quad \text{ for all }x,y\in D .
\end{equation}
To use Lemma \ref{l:estimate-for-halfspace-H3}, we further assume that \eqref{e:Hcond} holds when $d \ge 2$.
 Note that \eqref{e:Hcond} is clearly satisfied in case when $D$ is a half-space. 
We show now that under the above 
conditions, \textbf{(H2)}, \textbf{(H3)} and \textbf{(H5)} also hold.

\medskip
Let $a\in (0,1)$ and $x,y\in D$ such that $\delta_D(x)\wedge \delta_D(y)\ge a|x-y|$. 
Then 
$|x-y|^2/\delta_D(x)\delta_D(y))\le  1/a^2 $.  Since $\Phi$ is bounded on $[0, 1/a^2  )$, \textbf{(H2)} holds true.

\medskip
Let $a\in (0,1/2]$ and $|x-y|>a\delta_D(x)$. 
Then $\delta_D(y) \le|x-y|+ \delta_D(x) \le ((a+1)/a) |x-y|$.
Hence $|x-y|^2\ge  (a^2/(1+a)) \delta_D(x) \delta_D(y)$.
Therefore by
\eqref{e:Phi-infty}  and the fact that $\Phi(t)  \equiv \Phi(2)>0$ on $[0, 2)$, we conclude that there exists 
$c_1=c_1(a)>0$ 
such that
$$
\Phi\left(\frac{|x-y|^2}{\delta_D(x) \delta_D(y)}\right)\le 
c_1 \frac{|x-y|^{2\beta_2}}{\delta_D(x)^{\beta_2}\delta_D(y)^{\beta_2}}.
$$
Thus \eqref{e:assumption-B-example-2} holds, so \textbf{(H3)} follows form Lemma \ref{l:estimate-for-halfspace-H3}.

\medskip
Let $\epsilon \in (0,1)$,  $x_0\in D$ and $r>0$ with $B(x_0, (1+\epsilon)r)\subset D$. For $x_1,x_2\in B(x_0,r)$ and $z\in D\setminus B(x_0, (1+\epsilon)r)$, it holds that $|x_1-z|\le (1+2\epsilon)|x_2-z|$ and $\delta_D(x_1)\le \delta_D(x_2)+2r\le \delta_D(x_1)+(2/\epsilon)
\delta_D(x_2)=(1+2\epsilon)\delta_D(x_2)$. Therefore, there exists 
$c_2=c_2(\epsilon)\ge 1$ such that
$$
c_2^{-1}\frac{|x_1-z|^2}{\delta_D(x_1) \delta_D(z)}\le \frac{|x_2-z|^2}{\delta_D(x_2) \delta_D(z)} \le c_2 \frac{|x_1-z|^2}{\delta_D(x_1) \delta_D(z)}.
$$
Using this, \eqref{e:Phi-infty}  and the fact that $\Phi(t)  \equiv \Phi(2)>0$ on $[0, 2)$, \textbf{(H5)} holds true.

\vspace{.1in}
\textbf{Acknowledgment}:
We thank the referee for helpful comments and suggestions.

\vskip 0.1truein

\parindent=0em

{\bf Panki Kim}

Department of Mathematical Sciences and Research Institute of Mathematics,

Seoul National University, Seoul 08826, Republic of Korea

E-mail: \texttt{pkim@snu.ac.kr}

\bigskip

{\bf Renming Song}

Department of Mathematics, University of Illinois, Urbana, IL 61801,
USA

E-mail: \texttt{rsong@illinois.edu}

\bigskip

{\bf Zoran Vondra\v{c}ek}

Department of Mathematics, Faculty of Science, University of Zagreb, Zagreb, Croatia
Email: \texttt{vondra@math.hr}

\end{document}